\documentclass[12pt,reqno]{amsart}
\usepackage{amssymb,amscd,amsbsy}
\newcommand{\lb}{\linebreak}

\renewcommand{\a}{\alpha}
\renewcommand{\b}{\beta}

\newcommand{\e}{\varepsilon}
\newcommand{\vk}{\varkappa}
\newcommand{\z}{\zeta}

\newcommand{\vt}{\vartheta}
\renewcommand{\l}{\lambda}
\newcommand{\ro}{\rho}

\renewcommand{\t}{\tau}

\newcommand{\f}{\varphi}
\renewcommand{\o}{\omega}

\newcommand{\D}{\Delta}
\renewcommand{\L}{\Lambda}

\renewcommand{\O}{\Omega}
\newcommand{\U}{\Upsilon}

\newcommand{\A}{{\mathcal A}}
\newcommand{\B}{{\mathcal B}}

\newcommand{\F}{{\mathcal F}}
\newcommand{\Q}{{\mathcal Q}}
\newcommand{\h}{{\mathcal H}}

\newcommand{\cL}{{\mathcal L}}

\newcommand{\cR}{{\mathcal R}}

\newcommand{\cU}{{\mathcal U}}
\newcommand{\V}{{\mathcal V}}
\newcommand{\W}{{\mathcal W}}

\newcommand{\1}{{\bf 1}}

\newcommand{\C}{{\Bbb C}}
\newcommand{\T}{{\Bbb T}}
\newcommand{\pp}{{\Bbb P}}
\newcommand{\dd}{{\Bbb D}}

\newcommand{\mm}{{\Bbb M}}
\newcommand{\0}{{\boldsymbol{0}}}

\newcommand{\bs}{\boldsymbol}

\newcommand{\m}{{\boldsymbol m}}
\newcommand{\bS}{{\boldsymbol S}}

\newcommand{\rf}[1]{(\ref{#1})}

\newcommand{\df}{\stackrel{\mathrm{def}}{=}}
\newcommand{\dist}{\operatorname{dist}}

\newcommand{\trace}{\operatorname{trace}}
\newcommand{\rank}{\operatorname{rank}}

\newcommand{\eeq}{\end{equation}}
\newcommand{\beq}{\begin{equation}}
\newcommand{\bay}{\begin{eqnarray}}
\newcommand{\ba}{\begin{align*}}
\newcommand{\ea}{\end{align*}}
\newcommand{\ey}{\end{eqnarray}}
\newcommand{\bey}{\begin{eqnarray*}}
\newcommand{\eey}{\end{eqnarray*}}

\newcommand{\imp}{\Rightarrow}
\newcommand{\be}{\infty}

\newcommand{\bl}{\blacksquare}
\newcommand{\ess}{\operatorname{ess}}

\newcommand{\Range}{\operatorname{Range}}

\newcommand{\Pf}{{\bf Proof. }}

\newcommand{\ov}{\overline}

\newtheorem{thm}{\hspace{\parindent}Theorem}[section]

\newtheorem{cor}[thm]{\hspace{\parindent}Corollary}
\newtheorem{lem}[thm]{\hspace{\parindent}Lemma}

\begin{document}

\newcommand{\vse}{\vspace{.2in}}
\numberwithin{equation}{section}
\newcommand{\cG}{{\mathcal G}}
\newcommand{\cO}{{\mathcal O}}

\title{\bf Analytic approximation of matrix functions in $\bs{L^p}$}
\author{L. Baratchart, F.L. Nazarov, and V.V. Peller}
\thanks{The second author is partially supported by NSF grant DMS 0501067, the
third author is partially supported by NSF grant DMS 0700995}
\maketitle

\begin{abstract}
We consider the problem of approximation of matrix functions of class $L^p$ on the unit circle 
by matrix functions analytic in the unit disk in the norm of $L^p$, $2\le p<\be$. For an $m\times n$ matrix function $\Phi$ in $L^p$, we consider the Hankel operator $H_\Phi:H^q(\C^n)\to H^2_-(\C^m)$, $1/p+1/q=1/2$.
It turns out that the space of $m\times n$ matrix functions in $L^p$ splits into two subclasses: the set of respectable matrix functions and the set of weird matrix functions. If $\Phi$ is respectable, then its distance to the set of analytic matrix functions is equal to the norm of $H_\Phi$. For weird matrix functions, to obtain the distance formula, we consider Hankel operators defined on spaces of matrix functions. We also describe the set of $p$-badly approximable matrix functions in terms of special factorizations and give a parametrization formula for all best analytic approximants in the norm of
$L^p$. Finally, we introduce the notion of $p$-superoptimal approximation and prove the uniqueness of
a $p$-superoptimal approximant for rational matrix functions.
\end{abstract}

\section{\bf Introduction}
\setcounter{equation}{0}

\

The classical problem of analytic approximation of functions on the unit circle $\T$ is for a given 
function $\f\in L^\be$, to find a best $H^\be$ approximant to $\f$, i.e., to find a bounded analytic function $\psi$ in the unit disk $\dd$ such that 
$$
\|\f-\psi\|_{L^\be(\T)}=\dist_{L^\be}(\f,H^\be).
$$
A standard compactness argument shows that such a best approximant always exists, though it is not necessarily unique in general. However, under certain mild assumptions the best approximation is indeed unique. For example, this happens if $\f$ is continuous which was proved for the first time in \cite{Kh}.
We refer the reader to \cite{P} for a comprehensive study of the problem of best uniform approximation
by analytic functions.

It turns out that this approximation problem is closely related to Hankel operators
on the Hardy class $H^2$. For a function $\f\in L^\be$
the {\it Hankel operator} 
$$
H_\f:H^2\to H^2_-\df L^2\ominus H^2
$$
is defined by
$$
H_\f f=\pp_-\f f,\quad f\in H^2,
$$
where $\pp_-$ is the orthogonal projection from $L^2$ onto $H^2_-$. It was proved by Nehari 
(see \cite{P}, Ch. 1, \S\,1) that 
$$
\|H_\f\|=\dist_{L^\be}(\f,H^\be).
$$
Moreover, it turns out that the Hankel operators provide a powerful tool to constructively study the problem of
best uniform analytic approximation, see Chapters 1, 5, and 7 of \cite{P}. The problem of
uniform approximation by analytic functions is also called the {\it Nehari problem}.

The Nehari problem is very important in applications in control theory (see \cite{F} and \cite{P})
and is also a useful tool in
identification, see \cite{Pa} and \cite{BLPT}.
Moreover, for the needs of control theory it is important to consider not only the scalar case,
but also the case of matrix-valued functions. 

Let $\Phi$ be a bounded function with values in the space $\mm_{m,n}$ of $m\times n$
matrices (notationally, $\Phi\in L^\be(\mm_{m,n})$). The problem of best analytic approximation is to find 
a bounded analytic matrix function $Q$ of size $m\times n$ such that
$$
\|\Phi-Q\|_{L^\be(\mm_{m,n})}=\dist_{L^\be}\big(\Phi,H^\be(\mm_{m,n})\big),
$$
where $H^\be(\mm_{m,n})$ is the space of bounded analytic $m\times n$ matrix functions
and for a matrix function $\Psi\in L^\be(\mm_{m,n})$ we use the notation
$$
\|\Psi\|_{L^\be}\df\ess\sup_{\!\!\!\!\!\!\!\!\!\!\!\!\!\!\!z\in\T}\|\Psi(\z)\|_{\mm_{m,n}},
$$
where for a matrix $A$ in $\mm_{m,n}$ we denote by $\|A\|_{\mm_{m,n}}$
the operator norm of $A$ as an operator from $\C^n$ to $\C^m$.

As in the scalar case, the following distance formula holds:
$$
\dist_{L^\be}\big(\Phi,H^\be(\mm_{m,n})\big)=\|H_\Phi\|, \quad\Phi\in L^\be(\mm_{m,n}),
$$
where the Hankel operator $H_\Phi:H^2(\C^n)\to H^2_-(\C^m)\df L^2(\C^m)\ominus H^2(\C^m)$ is
defined by
$$
H_\Phi f=\pp_-\Phi f,\quad f\in H^2(\C^n),
$$
and $\pp_-$ is the orthogonal projection onto $H^2_-(\C^m)$ (see, e.g., \cite{P}, Ch. 2).

However, unlike the scalar case, even if $\Phi$ is a polynomial matrix function, generically $\Phi$
has infinitely many best approximants. To choose among all best approximants the ``very best approximant'', it is natural to consider the notion of superoptimal approximation. We refer the reader to 
\S\,2 of this paper for the definition of superoptimal approximation.

In this paper we are going to consider the problem of analytic approximation in the $L^p$ norm,
$2\le p<\be$.

Let $\f$ be a scalar function in $L^p$. The problem of best analytic approximation is to find a function $\psi$ in the Hardy class $H^p$ such that 
$$
\|\f-\psi\|_{L^p}=\dist_{L^p}(\f,H^p).
$$
If $1<p<\be$, then the space $L^p$ is uniformly convex which implies that every function 
$\f\in L^p$ has a unique best analytic approximant $\psi$ in the $L^p$. The function $\psi$ is said to be 
the {\it $p$-best analytic approximant to} $\f$. 

In \cite{BS} Hankel operators have been used to study the problem of best analytic 
and meromorphic approximation in $L^p$ for $2\le p<\be$ (see also \cite{Pr} for a dual approach). For $\f\in L^p$ the Hankel operator
$$
H_\f:H^q\to H^2_-
$$ 
is defined by 
$$
H_\f f=\pp_-\f f,\quad f\in H^q,
$$
where the exponent $q$ satisfies the equality
\bay
\label{pq2}
\frac1p+\frac1q=\frac12.
\ey
{\it Throughout this paper we always assume that $2\le p<\be$ and $q$ satisfies} \rf{pq2}.
{\it In the proofs of the results we assume that $2<p<\be$, it is an elementary exercise to adjust the
proofs for $p=2$.}

As in the case of uniform analytic approximation, the following formula holds
$$
\|H_\f\|_{H^q\to H^2_-}=\dist_{L^p}(\f,H^p)
$$
In \S\,2 of this paper we discuss in more detail the problem of best analytic approximation 
by scalar analytic functions in $L^p$.

In this paper we deal with the problem of approximation in $L^p$ by analytic matrix functions: given a function
$\Phi$ in $L^p(\mm_{m,n})$ (i.e., all entries of $\Phi$ belong to $L^p$), we search for a best analytic approximant $Q\in H^p(\mm_{m,n})$, i.e.,
$$
\|\Phi-Q\|_{L^p}=\dist_{L^p}\big(\Phi,H^p(\mm_{m,n})\big),
$$
where for a matrix function $\Psi\in L^p(\mm_{m,n})$,
$$
\|\Psi\|_{L^p}\df\|\Psi\|_{L^p(\mm_{m,n})}=\left(\int_\T\|\Psi(\z)\|_{\mm_{m,n}}^p\,d\m(\z)\right)^{1/p}.
$$
If we consider the Hankel operator
$$
H_\Phi:H^q(\C^n)\to H^2_-(\C^m)
$$
defined by
$$
H_\Phi f=\pp_-\Phi f,\quad f\in H^q(\C^n),
$$
it is easy to verify that 
$$
\|H_\Phi\|\le\dist_{L^p}\big(\Phi,H^p(\mm_{m,n})\big)
$$
(see Lemma \ref{inty}). It will be shown in \S\,4 that if $\Phi$ has a $p$-best analytic approximant
$Q$ such that for $\z$ in a subset of $\T$ of positive measure,
the space of maximizing vectors of $(\Phi-Q)(\z)$ is  one-dimensional,
then
\bay
\label{fla}
\|H_\Phi\|=\dist_{L^p}\big(\Phi,H^p(\mm_{m,n})\big).
\ey
Clearly, generically for an $m\times n$ matrix $A$, the maximizing vectors of $A$ span a one-dimensional subspace.

This makes it plausible that for a dense subset of matrix functions $\Phi$ in $L^p(\mm_{m,n})$
the distance formula \rf{fla} holds which would imply that \rf{fla} holds for all
 matrix functions $\Phi$ in $L^p(\mm_{m,n})$.
 
{\it Surprisingly, this is false!} 

In \S\,3 of this paper we obtain certain factorization theorems for analytic matrix functions that will be used to study Hankel operators. The main tool used in \S\,3 is Sarason's factorization theorem \cite{Sa}.

In \S\,4 we study the class of matrix functions $\Phi\in L^\be(\mm_{m,n})$, for which the distance formula \rf{fla} holds. Such matrix functions are called {\it respectable}. We obtain several characterizations of the class of respectable matrix functions.

The main result of \S\,5 is a construction of a $2\times2$ matrix function $\Phi$, for which 
\rf{fla} is false. Such matrix functions are called {\it weird}.

Thus the space $L^p(\mm_{m,n})$ splits in two subsets: the set of respectable 
matrix functions and the set of weird matrix functions. To compute the distance from 
a respectable matrix function $\Phi$ to the set of analytic matrix functions, we can use the distance formula \rf{fla}. However, to compute $\dist_{L^p}\big(\Phi,H^p(\mm_{m,n})\big)$ for weird matrix functions $\Phi$, we have to search for another formula. Note that in a sense both the set of respectable
matrix functions and the set of weird matrix functions are massive subsets of $L^p(\mm_{m,n})$; see the
discussion at the end of \S\,5.

It turns out, however, that the distance $\dist_{L^p}\big(\Phi,H^p(\mm_{m,n})\big)$ from $\Phi$ to the 
set of analytic matrix functions can be obtained 
for all matrix functions in $L^p$ as the norm of a Hankel operator if we consider Hankel operators acting on spaces of matrix functions rather than vector functions. Indeed, If we consider the Hankel operator $\bs{H}_\Phi$ defined on the space $H^q(\bS_2^n)$ of $n\times n$ matrix functions
with the norm
$$
\|F\|_{L^q(\bS_2^n)}=\left(\int_\T\|F(\z)\|_{\bS_2^n}^q\,d\m(\z)\right)^{1/q},
$$
Then the norm of the Hankel operator 
$$
\bs{H}_\Phi:H^q(\bS_2^n)\to H^2_-(\bS_2^n)
$$
is equal to $\dist_{L^p}\big(\Phi,H^p(\mm_{m,n})\big)$. Here for an $n\times k$ matrix $A$
the norm $\|A\|_{\bS_2^{n,k}}$ is the Hilbert--Schmidt norm of $A$ and
$\|A\|_{\bS_2^{n}}\df\|A\|_{\bS_2^{n,n}}$. This will be proved in \S\,6.
We also consider in $\S\,6$ Hankel operators acting on spaces of $n\times k$ matrix functions
and we introduce in \S\,6 the class of $n\times n$ matrix functions in $L^p$ of order $k$, $1\le k\le n$.

In \S\,7 we obtain a description of the set of $p$-badly approximable matrix functions. A matrix function $\Phi\in L^p(\mm_{m,n})$ is called {\it $p$-badly approximable} if 
$$
\|\Phi\|_{L^p}=\dist_{L^p}\big(\Phi,H^p(\mm_{m,n})\big).
$$
To obtain such a description, we use special factorizations that involve balanced matrix functions
(see \cite{P}, Ch. 14, \S\,1).

We also obtain in \S\,7 a parametrization formula for all $p$-best approximants.

In the last section we define the notion of $p$-superoptimal approximation and prove for rational matrix functions the uniqueness of a $p$-superoptimal approximant.

In \S\,2 we collect necessary information. In \S\,2.1 we present results on analytic approximation 
in $L^p$ of scalar functions. In \S\,2.2 we define the notion of superoptimal approximation and state some uniqueness results and properties of superoptimal approximants. Finally, in \S\,2.3 we define
the notion of balanced matrix functions and state factorization formulas for badly approximable matrix function.

Note that it suffices to study the problem of analytic approximation only for square matrix functions. Indeed, if a matrix function $\Phi$ is not square, we can add to $\Phi$ zero columns or zero rows to make it square. For the sake of simplicity, {\it beginning {\em\S\,6}, we state all the results only
for square matrix functions}.

\medskip

{\bf Notation and terminology.} Throughout the paper we are going to use the following notation and terminology:
\medskip

if $X$ and $Y$ are normed spaces and $T:X\to Y$ is a bounded linear operator, a vector $x\in X$ is called 
a {\it maximizing vector of} $T$ if 
$$
x\ne\0\quad\mbox{and}\quad\|Tx\|_Y=\|T\|\cdot\|x\|_X;
$$

if both $X$ and $Y$ are Hilbert spaces and $T$ is a bounded linear operator from $X$ to $Y$, then,
by definition, the {\it space of maximizing vectors of} $T$ is
$$
\{x\in X:~\|Tx\|_Y=\|T\|\cdot\|x\|_X\}
$$
(it is well known that the space of maximizing vectors is a closed subspace of $X$ that consists of the maximizing vectors and 
the zero vector);

$\mm_{m,n}$ is the space of $m\times n$ matrices;

$\mm_{n}\df\mm_{n,n}$;

if $X$ is a normed space of functions on $\T$, then $X(\mm_{m,n})$ means the space of $m\times n$ matrix functions whose entries belong to $X$. If this does not lead to a confusion, we say that
$\Phi\in X$ for an $m\times n$ matrix function $\Phi$ if $\Phi\in X(\mm_{m,n})$;

if $X=L^s$, $1\le s\le\be$, and $\Phi\in X(\mm_{m,n})$, then
$$
\|\Phi\|_X\df\|\Phi\|_{X(\mm_{m,n})}\df\|\ro\|_X,\quad
\mbox{where}\quad\ro(\z)\df\|\Phi(\z)\|_{\mm_{m,n}},\quad\z\in\T;
$$

for $1\le p\le\be$, the space $H^p(\mm_{m,n})$ is the subspace of $L^p(\mm_{m,n})$ that consists of matrix functions with entries in $H^p$. By definition, 
$$
H_0^p(\mm_{m,n})=\left\{F\in L^p(\mm_{m,n}):~F(0)=0\right\};
$$

for an operator $A$ on Hilbert space (or for a matrix $A$), the {\it singular values} $s_j(A)$
are defined by
$$
s_j(A)=\inf\{\|A-K\|:~\rank K\le j\};
$$

the Schatten--von Neumann class $\bS_r$, $1\le r<\be$, consists of operators $A$ on Hilbert space
with finite norm
\bay
\label{Scha}
\|A\|_{\bS_r}=\left(\sum_{j\ge0}s^r_j(A)\right)^{1/r};
\ey

for $r\in[1,\be)$, we denote by $\bS_r^{m,n}$ the space of $m\times n$ matrices $A$ equipped with the Schatten--von Neumann norm \rf{Scha};

$\bS_r^{n}\df\bS_r^{n,n}$;

if $X=L^s$, $1\le s\le\be$, then $X(\bS_r^{m,n})$ is the space of $m\times n$ matrix functions
with entries in $X$ equipped with the norm
$$
\|\Phi\|_{X(\bS_r^{m,n})}\df\|\ro\|_X,\quad\mbox{where}\quad
\ro(\z)\df\|\Phi(\z)\|_{\bS_r^{m,n}},\quad\z\in\T.
$$

\medskip

{\bf Acknowledgements.} We would like to thank Gilles Pisier and Ilya Spitkovskii for helpful remarks.
We are also grateful to Alberto Condori who has thoroughly read the manuscript and suggested several corrections.

\

\section{\bf Preliminaries}
\setcounter{equation}{0}

\

{\bf2.1. Analytic approximation by scalar functions in $\bs{L^p}$, $\bs{2\le p<\be}$.}  The problem of analytic approximation in $L^p$ was studied by many mathematicians, see, e.g., \cite{Sh} and
\cite{Ka}.
As we have already mentioned in the introduction, in \cite{BS} to study the problem of best analytic approximation in $L^p$, Hankel operators from $H^q$ to $H^2_-$ were used, where the exponent $q$ satisfies 
\rf{pq2} (see also \cite{Pr} in which a similar approach is used). The approach of \cite{BS} and \cite{Pr} is based on the analog of Nehari's theorem:
\bay
\label{NTp}
\|H_\f\|_{H^q\to H^2_-}=\dist_{L^p}(\f,H^p),\quad\f\in L^p.
\ey
Moreover, it can be shown that if $\f\in L^2$, then the Hankel operator $H_\f$ defined on the set of analytic polynomials by the formula
$$
H_\f f=\pp_-\f f
$$
extends to a bounded operator from $H^q$ to $H^2_-$ if and only if $\pp_-\f\in L^p$. This can be proved in exactly the same way as in the case of classical Hankel operators from $H^2$ to $H^2_-$ (see, e.g.,
\cite{P}, Ch. 1, \S\,1).  In particular, this implies that all bounded Hankel operators from $H^q$ to $H^2_-$
are compact, since the trigonometric polynomials are dense in $L^p$ and $H_\f$ has finite rank if $\f$ is a trigonometric polynomial.

A scalar function $\f\in L^p$ is called {\it $p$-badly approximable} if 
$$
\|\f-\psi\|_{L^p}\ge\|\f\|_{L^p}
$$
for any $\psi\in H^p$.

The following result describes the class of all $p$-badly approximable functions.

\begin{thm}
\label{spba}
Let $\f$ be a nonzero function in $L^p$. Then $\f$ is $p$-badly approximable if and only if there
exists an inner function $\vt$ and an outer function $h$ in $H^2$ such that
\bay
\label{spbaf}
\f=\bar z\bar\vt\frac{\bar h}{h^{2/q}}=\bar z\bar\vt\frac{\bar h}{h^{\frac{p-2}{p}}}.
\ey
\end{thm}

\Pf Suppose that $\f$ is $p$-badly approximable.
Let $f\in H^q$ be a maximizing vector of $H_\f$. Such a vector exists, since $H_\f$ is compact.
We have
\begin{align*}
\|H_\f f\|_{L^2}&=\|\pp_-\f f\|_{L^2}\le \|\f f\|_{L^2}\\[.2cm]
&\le\|\f\|_{L^p}\|f\|_{L^q}
=\|H_\f\|\cdot\|f\|_{L^q}=\|H_\f f\|_{L^2},
\end{align*}
since $f$ is a maximizing vector. Thus all inequalities in the above chain are equalities.
The fact that the first inequality turns into equality means that $\f f\in H^2_-$. The second inequality turns into equality if and only if $|\f|^p=c|f|^q$ for some $c>0$. We can multiply $f$ by a constant after which 
$c$ becomes equal to 1. Let $h$ be an outer function in $H^2$ such that
$|h|=|f|^{q/2}$. Then $f$ admits a factorization $f=\vt_1h^{2/q}$, where $\vt_1$ is an inner function.

Put $g=\bar z\ov{H_\f f}\in H^2$. We have $|g|^2=|\f f|^2=|h|^2$. Let $g=\vt_2h$, where $\vt_2$
is an inner function. Then 
$$
\f=\frac{\bar z \bar g}{f}=\bar z\bar\vt_1\bar\vt_2\frac{\bar h}{h^{2/q}}.
$$
It remains to put $\vt=\vt_1\vt_2$.

Suppose now that $\f$ is of the form \rf{spbaf}. Put $f=h^{2/q}$. We have
$$
\|H_\f\|\cdot\|f\|_{L^q}\ge\|H_\f f\|_{L^2}=\|h\|_{L^2}=\|\f\|_{L^p}\|f\|_{L^q}\ge\|H_\f\|\cdot\|f\|_{L^q}.
$$
Thus $\|\f\|_{L^p}=\|H_\f\|$, and so $\f$ is $p$-badly approximable. $\bl$

\medskip

{\bf Remark.} Note that in the case $p=\infty$ the situation is slightly different. 
A bounded Hankel operators from $H^2$ to $H^2_-$ is not necessarily compact
and does not necessarily have a maximizing vector. A badly approximable function $\f$ has the form 
$$
\f=c\bar z\bar\vt\frac{\bar h}{h},
$$
where $c\in\C$, $\vt$ is an inner function, and $h$ is an outer function in $H^2$,
if and only if the Hankel operator $H_\f:H^2\to H^2_-$ has a maximizing vector, see \cite{P}, Ch. 1, \S\,1.

In the case $p=2$, Theorem \ref{spba} means that the $2$-badly approximable functions are precisely the functions in $H^2_-$ and
a function $f\in H^\be$ is a maximizing vector of the Hankel operator $H_\f:H^\be\to H^2_-$
with a $2$-badly approximable symbol $\f$ if and only if $f=c\vt$, where $c$ is a nonzero complex number and $\vt$ is an inner divisor of $\bar z\bar\f$.

\begin{cor}
\label{mod}
Let $\o$ be a nonnegative function in $L^p$. The following are equivalent:

{\em(i)} there exists a $p$-badly approximable function $\f\in L^p$ such that $|\f|=\o$;

{\em(ii)} $\log\o\in L^1$.
\end{cor}

\Pf The implication (i)$\imp$(ii) is an immediate consequence of Theorem \ref{spba} and the fact that the logarithm of the modulus
of any outer function is in $L^1$.

Conversely, suppose that $\log\o\in L^1$. Let $h$ be an outer function such that \lb$|h|=\o^{p/2}$. Clearly, $h\in H^2$.
Let $\f=\bar z\frac{\bar h}{h^{2/q}}$. By Theorem \ref{spba}, $\f$ is badly approximable. We have
$$
|\f|=|h|^{1-2/q}=|h|^{2/p}=\o.\quad\bl
$$

\begin{cor}
\label{desc}
Let $\o$ be a nonnegative function in $L^p$ such that $\log\o\in L^1$ and let $h$ be an outer function such that
$|h|=\o^{p/2}$. Then the set of $p$-badly approximable functions with modulus $\o$ coincides with
$$
\left\{\bar z\bar\vt\frac{\bar h}{h^{2/q}}:~\vt~\mbox{ is an inner function}\right\}.
$$
\end{cor}

\Pf The result follows immediately from Theorem \ref{spba}. $\bl$

\medskip

{\bf 2.2. Superoptimal approximation.} As we have  already mentioned in the introduction, even for polynomial matrix functions $\Phi$ there can be many best \lb analytic approximants in the $L^\be$ norm.
For instance, if $\Phi=\left(\begin{matrix}\bar z&\0\\\0&\0\end{matrix}\right)$ and 
\lb$F=\left(\begin{matrix}\0&\0\\\0&f\end{matrix}\right)$, where $f$ is a scalar function in the unit ball of 
$H^\be$, then $F$ is a best approximant to $\Phi$.

To introduce the notion of superoptimal approximation, recall the notion of singular values of matrices.
For a matrix $A$ the $j$th singular value of $A$
is defined by
$$
s_j(A)\df\inf\{\|A-K\|:~\rank K\le j\},\quad j\ge0.
$$
Clearly, $s_0(A)=\|A\|$.

\medskip

{\bf Definition.}
Given a matrix function $\Phi\in L^\be(\mm_{m,n})$ we define inductively
the sets $\bs{\O}_j$, $0\le j\le\min\{m,n\}-1$, by
$$
\bs{\O}_0=\left\{Q\in H^\be(\mm_{m,n})
:~Q~\mbox{minimizes}~\ t_0\df\ess\sup_{\z\in\T}\|\Phi(\z)-Q(\z)\|\right\};
$$
$$
\bs{\O}_j=\left\{Q\in \O_{j-1}:~Q~\mbox{minimizes}~\ 
t_j\df\ess\sup_{\z\in\T}s_j\big(\Phi(\z)-Q(\z)\big)\right\},\quad j>0.
$$
Functions in $\bigcap\limits_{k\ge 0} \bs{\O}_k = \bs{\O}_{\min\{m,n\}-1} $ are called {\it
superoptimal approximants} to $\Phi$ by bounded analytic matrix functions.  
The numbers
$t_j=t_j(\Phi)$ are called the {\it superoptimal singular values} of $\Phi$.
Note that the matrix functions in $\bs{\O}_0$ are just the best approximants by analytic matrix functions.

In other words, a superoptimal approximant minimizes the essential suprema of 
the singular values of $(\Phi-Q)(\z)$ lexicographically.

It was proved in \cite{PY} that if $\Phi\in(H^\be+C)(\mm_{m,n})$ (i.e., each entry of $\Phi$ is a sum of 
a a continuous function and an $H^\be$ function), then $\Phi$ has a unique superoptimal approximant.
Moreover, if $Q$ is the unique superoptimal approximant to $\Phi$, then
$$
s_j\big((\Phi-Q)(\z)\big)=t_j,\quad\z\in\T.
$$
Later in \cite{PT} the same results were obtained under a less restrictive assumption on $\Phi$.
We refer the reader to \cite{P}, Ch.~14 for a detailed presentation of the theory of superoptimal approximation.

\medskip

{\bf 2.3. Balanced matrix functions and factorizations of badly approximable matrix functions.} A matrix function $\Phi$ in $L^\be(\mm_{m,n})$ is called {\it badly approximable} if 
$$
\|\Phi\|_{L^\be}\le\|\Phi-Q\|_{L^\be}
$$
for any $Q\in H^\be(\mm_{m,n})$. 

A matrix function $\Phi$ is called {\it very badly approximable} if the zero matrix function is a superoptimal approximant to $\Phi$.

In \cite{PY} and \cite{AP} the set of badly approximable matrix
functions of class \lb$(H^\be+C)(\mm_{m,n})$ was described in terms of certain special factorizations
(see also \cite{PT2} in which a geometric description of very badly approximable matrix functions
was obtained).
Such factorizations involve certain special unitary-valued matrix functions (balanced matrix functions),
see \cite{P}, Ch. 14, \S\,1. To define balanced matrix functions, we have to introduce several notions.

 A matrix function $\Theta\in H^\be(\mm_{n,k})$ is called {\it inner} if on the unit circle
$\Theta^*\Theta=\bs{I}_k$, where $\bs{I}_k$ is the matrix function identically equal to the 
identity matrix $I_k$.

A matrix function $F\in H^2(\mm_{m,n})$ is called {\it outer} if the set
$$
\{F f:~f~~\mbox{is a polynomial in}~~H^2(\C^n)\}
$$
is dense in $H^2(\C^m)$.

Finally, a matrix function  $F\in H^2(\mm_{m,n})$ is called {\it co-outer} if the transposed
function $F^{\rm t}$ is outer.

It is well known (see, e.g., \cite{N} or \cite{SF}) that if $\Psi$ is a matrix function of class $H^2$, then $\Phi$ admits an  {\it inner-outer factorization}
$$
\Phi=\Theta F,
$$
where $\Theta$ is an inner matrix function and $F$ is an outer matrix function.

Let $k<n$ and let $\U$ be an $n\times k$ inner and co-outer matrix function. It is well known (see \cite{P}, Ch. 14, \S\,1 and \cite{H}, Ch. 9) that there exists an inner and co-outer matrix function $\Theta$ of size
$n\times(n-k)$ such that the matrix function
\bay
\label{bmf}
\V=\left(\begin{matrix}\U&\ov{\Theta}\end{matrix}\right)
\ey
takes unitary values or, in other words, is {\it unitary-valued}. Matrix functions of the form \rf{bmf}
are called {\it balanced matrix functions}. If we want to specify that the analytic part of $\V$ has $k$ columns, we say that $\V$ is a {\it $k$-balanced matrix function}. In the case $k=1$, $k$-balanced matrix functions are also called {\it thematic matrix functions}. If $k=n$ by a $k$-balanced matrix function, we mean a matrix function of the form $\t\bs{I}_n$, where $\t$ is a complex number of modulus 1.

Balanced matrix functions have many interesting properties, see \cite{P}, Ch. 14, \S\,1.
They have been used to obtain a description of badly approximable matrix functions, to parametrize
the set of best analytic approximants, to characterize very badly approximable matrix functions,
to prove the uniqueness of superoptimal approximants, and to construct the superoptimal approximant (see, \cite{PY}, \cite{PT}, \cite{AP}, and \cite{P}, Ch. 14).

In particular, it was shown in \cite{PY} (see also \cite{P}, Ch.~14, \S\,2) that
 if $\Phi$ is \lb a matrix function in $L^\be(\mm_{m,n})$ such that the Hankel operator 
\lb$H_\Phi:H^2(\C^n)\to H^2_-(\C^m)$ has a maximizing vector, then $\Phi$ is badly approximable
if and only if $\Phi$ admits a factorization
$$
\Phi=\|H_\Phi\|
W^*\left(\begin{array}{cc}\bar z\bar\vt \bar h/h&\0\\[.2cm]\0&\Phi_\#\end{array}\right)V^*,
$$
where $V$ and $W^{\rm t}$ are thematic matrix functions, $\vt$ is a scalar inner function,
$h$ is a scalar outer function in $H^2$, and $\Phi_\#$ is a matrix function of size $(m-1)\times(n-1)$
such that $\|\Phi_\#(\z)\|_{\mm_{m-1,n-1}}\le1$ almost everywhere on $\T$.

Another characterization of badly approximable functions was obtained in \cite{AP}
(see also \cite{P}, Ch.~14, \S\,15). Let $\Phi\in(H^\be+C)(\mm_{m,n})$ and let $k$ be the
number of superoptimal singular values $t_j(\Phi)$ equal to $t_0(\Phi)$ (in other words, $k$ is the multiplicity of the superoptimal singular value $t_0(\Phi)$). Then $\Phi$ is badly approximable if and only if
$$
\Phi=\|H_\Phi\|
\W^*\left(\begin{array}{cc}U&\0\\[.2cm]\0&\Phi_\#\end{array}\right)\V^*,
$$
where $U$ is a $k\times k$ very badly approximable unitary-valued function of class \lb$H^\be+C$,
$\V$ and $\W$ are $k$-balanced matrix functions, and $\Phi_\#$ is a matrix function
in $(H^\be+C)(\mm_{m-k,n-k})$ such that $\|\Phi_\#(\z)\|_{\mm_{m-k,n-k}}\le1$ almost everywhere on 
$\T$ and $\|H_{\Phi_\#}\|<1$. Actually, the condition $\Phi\in(H^\be+C)(\mm_{m,n})$ can be relaxed
(see \cite{AP} and \cite{P}, Ch.~14, \S\,15).

\

\section{\bf Factorization of analytic matrix functions}
\setcounter{equation}{0}

\

In this section we obtain several factorization theorems for analytic matrix functions
that will be used to study Hankel operators.

We are going to use the following result by D. Sarason that is an analog of Riesz factorization:

\medskip

{\bf Sarason's Theorem [Sa].} {\it Let $\h$ be a separable Hilbert space and 
let $\Psi$ be an analytic integrable $\B(\h)$-valued function on $\T$.
Then there exist analytic square integrable functions $\Q$ and $\cR$
such that 
\bay
\label{sar}
\Psi=\Q\cR,\quad\cR^*\cR=\big(\Psi^*\Psi\big)^{1/2},\quad
\mbox{and}\quad\Q^*\Q=\cR\cR^*\quad\mbox{a.e. on}~~\T.
\ey}

\medskip

The following theorem can be deduced easily from Sarason's theorem.
Recall that $2\le p<\be$ and $q$ satisfies \rf{pq2}; as usual, $p'$ is the dual exponent: $1/p+1/p'=1$.

\begin{thm}
\label{st}
Let $\Psi\in H^{p'}(\mm_{n})$. Then there exist matrix functions
\lb$F\in H^q(\mm_{n})$ and $G\in H^2(\mm_{n})$
such that
$$
\Psi=FG\quad\mbox{and}\quad\|\Psi\|_{L^{p'}(\bS_1)}=
\|F\|_{L^q(\bS_2)}\|G\|_{L^2(\bS_2)}.
$$
\end{thm}

\medskip

{\bf Proof of Theorem \ref{st}.}  Clearly, we may assume that $\Psi$ is a nonzero function.
Suppose that $\Q$ and $\cR$ satisfy the requirements of Sarason's theorem.
Let $h$ be a scalar outer function such that
\bay
\label{vnesh}
|h(\z)|=\|\Psi(\z)\|_{\bS_1^n}^{1/2-p'/2},\quad\z\in\T.
\ey
Put
$$
F=h\Q\quad\mbox{and}\quad G=\frac1h\cR.
$$
By Sarason's Theorem,
\begin{align*}
\|F\|_{L^q(\bS_2^n)}^q&=\int_\T|h(\z)|^q\|\Q(\z)\|^q_{\bS_2^n}\,d\m(\z)\\[.2cm]
&=\int_\T\|\Psi(\z)\|_{\bS_1^n}^{(1/2-p'/2)q+q/2}\,d\m(\z)\\[.2cm]
&=\int_\T\|\Psi(\z)\|_{\bS_1^n}^{p'}\,d\m(\z).
\end{align*}
Similarly,
\begin{align*}
\|G\|_{L^2(\bS_2^n)}^2&=\int_\T|h(\z)|^{-2}\|\cR(\z)\|^2_{\bS_2^n}\,d\m(\z)\\[.2cm]
&=\int_\T\|\Psi(\z)\|_{\bS_1^n}^{p'-1+1}\,d\m(\z)\\[.2cm]
&=\int_\T\|\Psi(\z)\|_{\bS_1^n}^{p'}\,d\m(\z).
\end{align*}
It follows that 
$$
\qquad\qquad
\|F\|_{L^q(\bS_2^n)}\|G\|_{L^2(\bS_2^n)}=\|\Psi\|_{L^{p'}(\bS_1^n)}^{p'/q+p'/2}=\|\Psi\|_{L^{p'}(\bS_1^n)}.
\qquad\qquad\bl
$$

We need a version of Theorem \ref{st} in the case $\Psi(\z)$ has rank $k$ for $\z\in\T$.
The following result can be deduced from Sarason's theorem.

\begin{thm}
\label{nk12}
Let $1\le k\le n$ and let $\Psi$ be a function in $H^1(\mm_{n})$ such that
\bay
\label{posm}
\rank \Psi(\z)=k\quad \mbox{on a subset of $\T$ of positive measure}
\ey
Then there exist matrix functions $\F\in H^2(\mm_{n,k})$ and $\cG\in H^2(\mm_{k,n})$
such that
\bay
\label{FG}
\Psi=\F\cG\quad\mbox{and}\quad\|\Psi(\z)\|_{\bS_1^n}=\|\F(\z)\|_{\bS_2^{n,k}}\|\cG(\z)\|_{\bS_2^{k,n}},
\quad\z\in\T.
\ey
\end{thm}

\Pf Clearly, each minor of $\Psi$ belongs to the Hardy class $H^s$ for some $s>0$.
It follows now from the uniqueness theorem for Hardy classes that condition \rf{posm} is equivalent to the fact that $\rank\Psi(\z)=k$ almost everywhere on $\T$.

Let $\Q$ and $\cR$ be $n\times n$ matrix functions satisfying the requirements of Sarason's theorem.
Then
\bay
\label{rpro}
\|\Psi(\z)\|_{\bS_1^n}=\|\Q(\z)\|_{\bS_2^n}\|\cR(\z)\|_{\bS_2^n}.
\ey

We need the following elementary lemma whose proof is given here for completeness.

\begin{lem}
\label{AB}
If $A$ and $B$ are operators on Hilbert space, $\rank AB=k$,
and $\|AB\|_{\bS_1}=\|A\|_{\bS_2}\|B\|_{\bS_2}$, then $\rank A=\rank B=k$.
\end{lem}

Let us first complete the proof of Theorem \ref{nk12}.

By Lemma \ref{AB}, \rf{rpro} implies that 
\bay
\label{krank}
\rank \Q(\z)=k\quad\mbox{and}\quad\rank\cR(\z)=k\quad\mbox{for almost all}\quad\z\in\T.
\ey
Consider the inner-outer factorization of $\cR$:
$$
\cR=\U\cG,
$$
where $\U$ is an inner matrix function and $\cG$ is an outer matrix function. It follows from \rf{krank}
that $\U$ has size $n\times k$ and $\cG$ has size $k\times n$. We can define now the function $\F$ by
$\F=\Q\U$.
Since $\U$ takes isometric values almost everywhere on $\T$, it follows that
$$
\|\F(\z)\|_{\bS_2^{n,k}}=\|\Q(\z)\|_{\bS_2^n}\quad\mbox{and}\quad
\|\cG(\z)\|_{\bS_2^{k,n}}=\|\cR(\z)\|_{\bS_2^n},
$$
and so \rf{FG} holds. $\bl$

\medskip

{\bf Proof of Lemma \ref{AB}.} Clearly, if $\rank A<k$ or $\rank B<k$, then $\rank AB<k$. 
Suppose now that the conclusion of the lemma is false. Without loss of generality we may assume that
$\rank A>k$. Let $P$ be the orthogonal projection onto $\Range AB$. Then $AB=PAB$. Clearly,
$\|A\|_{\bS_2}^2=\|PA\|_{\bS_2}^2+\|(I-P)A\|_{\bS_2}^2$. Since $\rank P=k$ and $\rank A>k$, it follows
that $\|PA\|_{\bS_2}<\|A\|_{\bS_2}$. Thus
$$
\|AB\|_{\bS_1}=\|PAB\|_{\bS_1}\le\|PA\|_{\bS_2}\|B\|_{\bS_2}<\|A\|_{\bS_2}\|B\|_{\bS_2}
=\|AB\|_{\bS_1}
$$
and we get a contradiction. $\bl$

We need the following consequence of Theorem \ref{nk12}.

\begin{thm}
\label{nkq2}
Let $\Psi\in H^{p'}(\mm_{n})$ such that
$$
\rank\Psi(\z)=k,\quad\z\in\T.
$$
Then there exist matrix functions
$F\in H^q(\mm_{n,k})$ and $G\in H^2(\mm_{k,n})$
such that
$$
\Psi=FG\quad\mbox{and}\quad\|\Psi\|_{L^{p'}(\bS_1^n)}=
\|F\|_{L^q(\bS_2^{n,k})}\|G\|_{L^2(\bS_2^{k,n})}.
$$
\end{thm}

\Pf As in the proof of Theorem \ref{st}, we put
$$
F=h\F,\quad G=\frac1{h}\cG,
$$
where $\F$ and $\cG$ are matrix functions satisfying the requirements of Theorem \ref{nk12}
and $h$ is a scalar outer function satisfying \rf{vnesh}. The fact that $F$ and $G$ satisfy 
the conclusions of the theorem is exactly the same as in the proof of Theorem \ref{nk12}.
$\bl$



The case of  matrix functions of rank 1 is of special interest. We treat this case separately, without
using Sarason's theorem.

\begin{lem}
\label{r1e}
Let $\Psi\in H^{1}(\mm_{m,n})$ such that
\bay
\label{rank}
\rank \Psi(\z)=1\quad \mbox{on a subset of $\T$ of positive measure}
\ey
Then there exist vector functions $u\in H^2(\C^m)$,
and $v\in H^2(\C^n)$,  such that
\bay
\label{tra}
\Psi(\z)=u(\z)v^{\rm t}(\z),\quad\z\in\T.
\ey
and
\bay
\label{rav}
\|u(\z)\|_{\C^m}=\|v(\z)\|_{\C^n}=\|\Psi(\z)\|_{\mm_{m,n}}^{1/2}.
\ey
\end{lem}

\Pf Condition \rf{rank} means that each $2\times2$ minor of $\Psi$ vanishes on a set of positive 
measure. Since $\Psi\in H^1(\mm_{m,n})$, it follows that all $2\times2$ minors of $\Psi$
are identically equal to zero. Thus condition \rf{rank} implies that $\rank\Psi(\z)=1$ almost 
everywhere on $\T$.

Let $h$ be an outer function such that 
$$
|h(\z)|^2=\|\Psi(\z)\|_{L^1(\mm_{m,n})},\quad\z\in\T,
$$
and let $G=h^{-1}\Psi$. Clearly, $G\in H^2(\mm_{m,n})$.
Consider the columns of $G$. Let $\cL$ be the invariant subspace of 
multiplication by $z$ on $H^2(\C^m)$ spanned by the columns of $G$. By the Beurling--Lax
theorem (see \cite{N}), there exists an inner function $\U$ of size
$m\times k$ such that $\cL=\U H^2(\C^k)$. Since $\rank G(\z)=1$ almost everywhere, it follows that
$k=1$. Then there exist functions $v_1,v_2,\cdots,v_n$ such that the columns of the matrix function $G$ are $v_1\U,v_2\U,\cdots,v_n\U$. Let
$$
v=\left(\begin{array}{c}v_1\\v_2\\\vdots\\v_n\end{array}\right).
$$
Clearly, $v\in H^2(\C^n)$ and $G=\U v^{\rm t}$. It remains to put $u=h\U$ and observe that
$u\in H^2(\C^m)$ and both \rf{tra} and \rf{rav} hold. $\bl$

\begin{thm}
\label{ats}
Let $\Psi$ be a rank one matrix function in $H^{p'}(\mm_{m,n})$. Then there exist
column functions $f\in H^q(\C^m)$ and $g\in H^2(\C^n)$ such that 
\bay
\label{fgt}
\Psi=fg^{\rm t}\quad\mbox{and}\quad \|\Psi\|_{L^{p'}(\mm_{m,n})}=\|f\|_{L^q(\C^m)}\|g\|_{L^2(\C^n)}.
\ey
\end{thm}

\Pf Let $u$ and $v$ be the column functions satisfying \rf{tra} and \rf{rav}. Let $h$ be a scalar outer function satisfying \rf{vnesh}. Put
$$
f=hu\quad\mbox{and}\quad g=\frac1h v.
$$
It is easy to verify that $f\in H^q(\C^m)$, $g\in H^2(\C^n)$, and
the equalities in  \rf{fgt} hold. $\bl$

\

\section{\bf Respectable matrix functions}
\setcounter{equation}{0}

\

The main result of this section is Theorem \ref{Rf}, which gives us several characterizations
of the set of matrix functions $\Phi\in L^p(\mm_{m,n})$, for which 
$\dist_{L^p}\big(\Phi,H^p(\mm_{m,n})\big)$ is equal to the norm of the Hankel operator 
$H_\Phi:H^q(\C^n)\to H^2_-(\C^m)$. The description of this class of matrix functions
(such matrix functions will be called {\it respectable}) makes it very natural to hope that 
all matrix functions in $L^p(\mm_{m,n})$ are respectable. However, it will be shown in \S\,5 that this is 
not true.

\medskip

{\bf Definition.} A matrix function $\Phi\in L^p(\mm_{m,n})\setminus H^p(\mm_{m,n})$ is called {\it regularly approximable} if there exists a best approximant $Q\in H^p(\mm_{m,n})$ such that the space of maximizing vectors of $(\Phi-Q)(\z)$
is one-dimensional on a subset of $\T$ of positive measure.

\medskip

It follows from the Hahn--Banach theorem that for $\Phi\in L^p(\mm_{m,n})$,
$$
\dist_{L^p}\big(\Phi,H^p(\mm_{m,n})\big)=
\sup\left|\int_\T \trace\big(\Phi(\z)\Psi(\z)\big)\,d\m(\z)\right|
$$
where the supremum is taken over all $\Psi\in H^{p'}_0(\mm_{n,m})$ (i.e., $\Psi\in H^{p'}(\mm_{n,m})$  and
$\Psi(0)=0$) such that $\|\Psi\|_{L^{p'}(\bS_1^{n,m})}\le1$.

Since the space $L^p(\mm_{m,n})$ is reflexive, it follows that for a matrix function \lb$\Phi\in L^p(\mm_{m,n})\setminus H^p(\mm_{m,n})$ there exists a matrix function 
$\Psi\in H^{p'}_0(\mm_{n,m})$ such that
\bay
\label{de}
\|\Psi\|_{L^{p'}(\bS_1^{n,m})}=1\quad \!\mbox{and}\quad\!
\int_\T \trace\big(\Phi(\z)\Psi(\z)\big)\,d\m(\z)=\dist_{L^p}\big(\Phi,H^p(\mm_{m,n})\big).
\ey
Such a function $\Psi$ is called a {\it dual extremal function} of $\Phi$.

Recall that for a matrix function $\Phi\in L^p(\mm_{m,n})$, we consider the Hankel operator
\lb$H_\Phi:H^q(\C^n)\to H^2_-(\C^m)$ defined by
$$
H_\Phi f=\pp_-\Phi f,\quad f\in H^q(\C^n),
$$
where $1/p+1/q=1/2$. 

As we have mentioned in \S\,2, for Hankel operators with scalar symbols, formula \rf{NTp} holds. Thus it is easy to see that the norm of the Hankel operator 
\lb$H_\Phi:H^q(\C^n)\to H^2_-(\C^m)$ is equivalent to the distance in $L^p$ from $\Phi$ to $H^p(\mm_{m,n})$.
Since in the case of scalar symbols all bounded Hankel operators from $H^q$ to $H^2_-$ are
compact, we can obtain the following result.

\begin{lem}
\label{comp}
For an arbitrary matrix function $\Phi$ in $L^p(\mm_{m,n})$, the Hankel operator
$H_\Phi:H^q(\C^n)\to H^2_-(\C^m)$ is compact.
\end{lem}

\begin{cor}
\label{max}
Let $\Phi\in L^p(\mm_{m,n})$. Then $H_\Phi$ has a maximizing vector in $H^q(\C^n)$.
\end{cor}

The following lemma gives us an upper estimate for the norm of $H_\Phi$.

\begin{lem} 
\label{inty}
Let $\Phi\in L^p(\mm_{m,n})$. Then
$$
\|H_\Phi\|\le\dist_{L^p}\big(\Phi,H^p(\mm_{m,n})\big).
$$
\end{lem}

\Pf Since $H_{\Phi-Q}=H_\Phi$ for an arbitrary $Q$ in $H^p(\mm_{m,n})$, it suffices to prove the
inequality
$$
\|H_\Phi\|\le\|\Phi\|_{L^p(\mm_{m,n})},\quad \Phi\in L^p(\mm_{m,n}).
$$
Suppose that $f\in H^q(\C^n)$ and $g\in H^2_-(\C^m)$. We have by H\"older's inequality,
\begin{align*}
|(H_\Phi f,g)|&=|(\Phi f,g)|\le\int_\T|\Phi fg^*|\,d\m\\[.2cm]
&\le\left(\int_\T\|\Phi(\z)\|^p_{\mm_{m,n}}\,d\m(\z)\right)^{1/p}
\left(\int_\T\|f(\z)\|_{\C^n}^q\right)^{1/q}
\left(\int_\T\|g(\z)\|_{\C^m}^2\right)^{1/2}\\[.2cm]
&=\|\Phi\|_{L^p(\mm_{m,n})}\|f\|_{L^q(\C^n)}\|g\|_{L^2(\C^n)}.\quad\bl
\end{align*}

The following theorem gives us several characterizations of the class of matrix functions $\Phi$, for which
$\|H_\Phi\|=\dist_{L^p}\big(\Phi,H^p(\mm_{m,n})\big)$.

\begin{thm}
\label{Rf}
Let $\Phi\in L^p(\mm_{m,n})\setminus H^p(\mm_{m,n})$. The following are equivalent:

{\em(i)} $\|H_\Phi\|=\dist_{L^p}\big(\Phi,H^p(\mm_{m,n})\big)$;

{\em(ii)} $\Phi$ belongs to the closure of the set of regularly approximable functions in $L^p$;

{\em(iii)} $ \Phi$ has a dual extremal function $\Psi$ such that $\rank\Psi(\z)=1$ on a set of positive measure;

{\em(iv)} $ \Phi$ has a dual extremal function $\Psi$ such that $\rank\Psi(\z)=1$, $\z\in\T$;

{\em(v)} if $Q$ is a best approximant to $\Phi$, then $\Phi-Q$ admits a factorization
\bay
\label{pT}
\Phi-Q=
W^*\left(\begin{array}{cc}\bar z\bar\vt \bar h/h^{2/q}&\0\\[.2cm]\0&\Phi_\#\end{array}\right)V^*,
\ey
where  $V$ and $W^{\rm t}$
are thematic matrix functions, $\vt$ is a scalar inner function, 
$h$ is a scalar outer function in $H^2$, and $\Phi_\#$ is an $(m-1)\times(n-1)$
matrix function such that $\|\Phi_\#(\z)\|_{\mm_{m-1,n-1}}\le |h(\z)|^{2/p}$, $\z\in\T$.
\end{thm}

Note that in \rf{pT} the outer function $h$ must satisfy the equality 
$$
|h(\z)|^{2/p}=\|(\Phi-Q)(\z)\|_{\mm_{m,n}},\quad
\z\in\T.
$$

\medskip

{\bf Remark.} Since the set of matrices, for which the space of maximizing vectors is one-dimensional is dense in the space of matrices, this suggests a hope that the set of regularly approximable $m\times n$ matrix functions is dense in $L^p(\mm_{m,n})$. If this were true, then the distance formula
$\|H_\Phi\|=\dist_{L^p}\big(\Phi,H^p(\mm_{m,n})\big)$ would hold for an arbitrary matrix functions in
$L^p(\mm_{m,n})$. {\it Surprisingly, we will show in \S\,5 that this is not the case.}

\medskip

{\bf Definition.}
Matrix functions $\Phi\in L^p(\mm_{m,n})\setminus H^p(\mm_{m,n})$ satisfying one of the conditions 
(i)--(v) in the statement of Theorem 
\ref{Rf} are called {\it respectable matrix functions}. If a matrix function 
$\Phi\in L^p(\mm_{m,n})\setminus H^p(\mm_{m,n})$ is not respectable,
it is called a {\it weird function}.

\medskip

It follows immediately from Theorem \ref{Rf} that the set of respectable functions is closed in $L^p$, while
the set of weird functions is open.

\medskip

{\bf Proof of Theorem \ref{Rf}.} We start with the proof of the implication
(iv)$\imp$(i). Let $\Psi$ be a dual extremal function such that 
$\rank\Psi(\z)=1$, $\z\in\T$. Then $\Psi$ satisfies \rf{de}.

Since $\|H_\Phi\|$ is always less than or equal to 
$\dist_{L^p}(\Phi,H^p(\mm_{m,n})\big)$, we have to show that
$$
\|H_\Phi\|\ge\dist_{L^p}(\Phi,H^p(\mm_{m,n})\big).
$$
By Theorem \ref{ats}, there exist functions $f\in H^2(\C^n)$ and $g\in H^2_0(\C^m)$ such that
$$
\Psi=fg^{\rm t}\quad\mbox{and}\quad 1=\|\Psi\|_{L^{p'}(\mm_{m,n})}=\|f\|_{L^q(\C^m)}\|g\|_{L^2(\C^n)}.
$$
Without loss of generality we may assume that $\|f\|_{L^q(\C^n)}=\|g\|_{L^2(\C^m)}=1$.
We have
\begin{align*}
\|H_\Phi\|&\ge|(H_\Phi f,g^*)|=\left|\int_\T\trace\big((H_\Phi f)g^{\rm t}\big)\,d\m\right|\\[.2cm]
&=\left|\int_\T\trace(\Phi fg^{\rm t})\,d\m\right|
=\int_\T \trace\big(\Phi\Psi\big)\,d\m\\[.2cm]
&=\dist_{L^p}\big(\Phi,H^p(\mm_{m,n})\big).
\end{align*}

\medskip

Next, let us show that (i)$\imp$(v). Let $f\in H^q(\C^n)$ be a maximizing vector of $H_\Phi$
and let $Q\in H^p(\mm_{m,n})$ be a best approximant to $\Phi$. We have 
\begin{align*}
\|H_\Phi f\|_{L^2(\C^m)}&=\|H_{\Phi-Q} f\|_{L^2(\C^m)}=\|\pp_-(\Phi-Q) f\|_{L^2(\C^m)}
\le\|(\Phi-Q) f\|_{L^2(\C^m)}\\[.2cm]
&\le\|(\Phi-Q)\|_{L^p(\mm_{m,n})}\|f\|_{H^q(\C^n)}=\|H_\Phi\|\cdot\|f\|_{H^q(\C^n)}=\|H_\Phi f\|_{L^2(\C^m)}.
\end{align*}
Hence, both inequalities are equalities. The fact that the first inequality turns into equality means that
$(\Phi-Q)f\in H^2_-(\C^m)$. The fact that the second inequality turns into equality means that
$f(\z)$ is a maximizing vector of $(\Phi-Q)(\z)$ for almost all $\z\in\T$ and 
\begin{align*}
&\int_\T\Big(\|(\Phi-Q)(\z)\|_{\mm_{m,n}}\|f(\z)\|_{\C^n}\Big)^2d\m(\z)\\=&
\left(\int_\T\|(\Phi-Q)(\z)\|_{\mm_{m,n}}^pd\m(\z)\right)^{2/p}
\left(\int_\T\|f(\z)\|_{\C^n}^qd\m(\z)\right)^{2/q},
\end{align*}
i.e., the corresponding H\"{o}lder inequality turns into equality, which implies that 
$\|(\Phi-Q)(\z)\|_{\mm_{m,n}}=c\|f(\z)\|^{q/p}_{\C^n}$ for some constant $c$. Since 
$$
\|(H_\Phi f)(\z)\|_{\C^m}=\|(\Phi-Q)(\z)\|_{\mm_{m,n}}\|f(\z)\|_{\C^n},
$$ 
it follows that $\|(H_\Phi f)(\z)\|_{\C^m}=c\|f(\z)\|^{q/2}_{\C^n}$. Multiplying the maximizing vector $f$ by a suitable constant,
one can always make the constant $c$ equal to 1, and so we may assume that
$$
\|(H_\Phi f)(\z)\|_{\C^m}=\|f(\z)\|^{q/2}_{\C^n}.
$$
Let $h$ be a scalar outer function such that 
$$
|h(\z)|=\|(H_\Phi f)(\z)\|_{\C^m},\quad\z\in\T,
$$
and so
$$
\|(\Phi-Q)(\z)\|_{\mm_{m,n}}=|h(\z)|^{2/p},\quad\z\in\T.
$$

 Put $g=\bar z\ov{H_\Phi f}\in H^2(\C^m)$. Then 
 $$
 \|g(\z)\|_{\C^m}=|h(\z)|\quad\mbox{and}\quad
 \|f(\z)\|_{\C^n}=|h(\z)|^{2/q}, \quad\z\in\T.
 $$
 The vector function $f$ admits 
a factorization $f=\vt_1h^{2/q}\bs{v}$, where $\vt_1$ is a scalar inner function and $\bs{v}$ is an $n\times1$ inner and co-outer function,
while the vector function $g$ admits a factorization 
$g=\vt_2h\bs{w}$, where $\vt_2$ is a scalar inner function and $\bs{w}$ is an $m\times1$ 
inner and co-outer function. 

Let now 
\bay
\label{tema}
V=\left(\begin{array}{cc}\bs{v}&\ov{\Theta}\end{array}\right)\quad \mbox{and}\quad 
W^{\rm t}=\left(\begin{array}{cc}\bs{w}&\ov{\Xi}\end{array}\right)
\ey
be thematic matrix functions
(see \S\,2.3).

Consider the matrix function $W(\Phi-Q)V$. Its upper left entry is equal to
\begin{align*}
\label{xi}
\xi&=\bs{w}^{\rm t}(\Phi-Q)\bs{v}=\bar\vt_2h^{-1}g^{\rm t}(\Phi-Q)\bar\vt_1h^{-2/q}f=
\bar\vt_1\bar\vt_2h^{-2/q-1}g^{\rm t}H_\Phi f\\[.2cm]
&=\bar z\bar\vt h^{-2/q-1}g^{\rm t}\bar g=\bar z\bar\vt h^{-2/q-1}|h|^2
=\bar z\bar\vt\frac{\bar h}{h^{2/q}}=\bar z\bar\vt\frac{\bar h}{h^{\frac{p-2}{p}}},
\end{align*}
where $\vt=\vt_1\vt_2$.

We have $|\xi(\z)|=\|(\Phi-Q)(\z)\|_{\mm_{m,n}}$. Since both $V$ and $W$ are unitary-valued,
it is easy to see that $\Phi-Q$ has the form \rf{pT}. 

\medskip

To prove the implication (v)$\imp$(ii), we need the following lemma. 

\begin{lem}
\label{lem}
Suppose that $\Phi$ is a matrix function that admits a factorization
\bay
\label{tfof}
\Phi=
W^*\left(\begin{array}{cc}\bar z\bar\vt \bar h/h^{2/q}&\0\\[.2cm]\0&\Phi_\#\end{array}\right)V^*,
\ey
where $\vt$, $h$, $\Phi_\#$, $V$, and $W$ are as in the statement of Theorem {\em\ref{Rf}}.
Then $\Phi$ is $p$-badly approximable.
\end{lem}

\Pf As we have already observed, for an arbitrary matrix function $\Phi$ in $L^p(\mm_{m,n})$ the following inequalities hold:
$$
\|H_\Phi\|\le\dist_{L^p}\big(\Phi,H^p(\mm_{m,n})\big)\le\|\Phi\|_{H^p(\mm_{m,n})}.
$$
It suffices to prove that if $\Phi$ is as in \rf{tfof}, then $\|H_\Phi\|\ge\|\Phi\|_{H^p(\mm_{m,n})}$. 
Consider the matrix functions $V$ and $W$:
$$
V=\left(\begin{matrix}\bs{v}&\Theta\end{matrix}\right)\quad\mbox{and}\quad
W^{\rm t}=\left(\begin{matrix}\bs{w}&\Xi\end{matrix}\right)
$$
Let $f=h^{2/q}\bs{v}$. It is easy to verify that 
$$
H_\Phi f=\Phi f=\bar z\bar h\bar\vt\ov{\bs{w}}\quad\mbox{and}\quad
\|H_\Phi f\|_{L^2(\C^m)}=\|\Phi\|_{L^p(\mm_{m,n})}\|f\|_{L^q(\C^n)}
$$
which implies that $\|H_\Phi\|\ge\|\Phi\|_{L^p(\mm_{m,n})}$. $\bl$

\medskip

(v)$\imp$(ii). Let $R=\Phi-Q$. For $\e>0$ we consider the function $R_\e$ defined by 
$$
R_\e=W^*\left(\begin{array}{cc}(1+\e)\bar z\bar\vt \bar h/h^{2/q}&\0\\[.2cm]\0&\Phi_\#\end{array}\right)V^*.
$$
By Lemma \ref{lem}, $R$ and $R_\e$ are $p$-badly approximable matrix functions. We define the 
function $\Phi_\e$ by $\Phi_\e=R_\e+Q$. 

Since $Q\in H^p(\mm_{m,n})$ and $R_\e$ is $p$-badly approximable, it follows that $Q$ is a
$p$-best approximant to $\Phi_\e$. Clearly, for $\z\in\T$, the space of maximizing vectors
of $R_\e(\z)$ is one-dimensional, and so $\Phi_\e$ is a regularly approximable matrix function.
The result follows from the obvious fact that 
$$
\|\Phi_\e-\Phi\|_{L^p}\to0\quad\mbox{as}\quad\e\to\be.
$$

\medskip

To show that (iii)$\imp$(iv), we observe that (iii) implies that each $2\times 2$ minor of $\Psi$
vanishes on a set of positive measure. By the uniqueness theorem for the Hardy classes, it follows that
all $2\times2$ minors of $\Psi$ are zero almost everywhere on $\T$ which proves (iv).

\medskip

Let us prove now that (ii)$\imp$(i). Clearly, it suffices to show that if $\Phi$ is regularly approximable, then $\|H_\Phi\|=\dist_{L^p}\big(\Phi,H^p(\mm_{m,n})\big)$. Let $Q$ be a matrix function in $H^p(\mm_{m,n})$ such that the space of maximizing vectors of $(\Phi-Q)(\z)$ is one-dimensional on 
a subset of $\T$ of positive measure. Let $\Psi$ be a dual extremal function of $\Phi$. It follows easily from \rf{de} that
\bay
\label{defn}
\trace\big((\Phi-Q)(\z)\Psi(\z)\big)=\|(\Phi-Q)(\z)\|_{\mm_{m,n}}\|\Psi(\z)\|_{\bS_1^{n,m}},\quad\z\in\T.
\ey
We need the following elementary lemma.

\begin{lem}
\label{elem}
Let $A\in\mm_{m,n}$ and $B\in\mm_{n,m}$ be matrices satisfying
$$
|\trace(AB)|=\|A\|_{\mm_{m,n}}\|B\|_{\bS_1^{n,m}}.
$$
Assume that the space of maximizing vectors $A$ is one-dimensional. Then
$B$ has rank $1$.
\end{lem}

\Pf Without loss of generality we may assume that $m=n$.
By considering the polar decomposition of $B$, we may assume 
that $B$  is positive, i.e., $(Bx,x)\ge0$ for every vector $x$.
Let $e_1,\cdots,e_n$ be an orthonormal basis of eigenvectors of $B$ and let $Be_j=\l_je_j$.
We have
\begin{align*}
|\trace(AB)|&=|\trace(BA)|=\left|\sum_{j=1}^n(BAe_j,e_j)\right|=\left|\sum_{j=1}^n(Ae_j,Be_j)\right|\\[.2cm]
&=\left|\sum_{j=1}^n\l_j(Ae_j,e_j)\right|\le\sum_{j=1}^n\l_j\|Ae_j\|\le\|A\|\sum_{j=1}^n\l_j.
\end{align*}
On the other hand,
$$
\|A\|\cdot\|B\|_{\bS_1^{n}}
=\|A\|\cdot\sum_{j=1}^n(Be_j,e_j)=\|A\|\sum_{j=1}^n\l_j.
$$
It follows that if $\|Ae_j\|<\|A\|$, then $\l_j=0$. By the hypotheses there can be only one $j$, for which 
$\|Ae_j\|=\|A\|$, which proves the result. $\bl$

\medskip

It follows from \rf{defn} and from Lemma \ref{elem} that $\Phi$ satisfies (iii). Since we have already proved that (iii)$\imp$(iv) and (iv)$\imp$(i), it follows that $\Phi$ satisfies (i).

\medskip

The fact that (iv)$\imp$(iii) is obvious. It remains to prove that (v)$\imp$(iv).
Suppose that $\Phi-Q$ is factorized as in \rf{pT}. Without loss of generality we may assume that
$\|\Phi-Q\|_{L^p}=1$
Define the matrix function $\Psi$ by
$$
\Psi=z\vt h^{1+2/q}\left(\begin{matrix}\bs{v}&\0\end{matrix}\right)
\left(\begin{matrix}\bs{w}^{\rm t}\\\0\end{matrix}\right),
$$
where $\bs{v}$ and $\bs{w}$ are as in \rf{tema}.

Clearly, $\rank\Psi(\z)=1$, $\z\in\T$. We have
$$
\|\Psi\|_{L^{p'}(\bS_1^{n,m})}^{p'}=\int_\T|h(\z)|^{p'(1+2/q)}\,d\m=\|h\|^2_{L^2}=1
$$
and
\begin{align*}
\int_\T\trace\big((\Phi-Q)\Psi\big)\,d\m&=
\int_\T z\vt h^{1+2/q}\trace\left(\left(\begin{matrix}\bs{w}^{\rm t}
\\[.2cm]\0\end{matrix}\right)(\Phi-Q)
\left(\begin{matrix}\bs{v}&\0\end{matrix}\right)\right)\,d\m\\[.2cm]
&=\int_\T\trace\left(\begin{matrix}|h|^2&\0\\[.2cm]\0&\0\end{matrix}\right)
\,d\m=\|h\|_{L^2}^2=\|\Phi-Q\|_{L^p}=1.
\end{align*}
This completes the proof. $\bl$

\medskip

{\bf Remark.} Note that in the case of analytic matrix approximation in the $L^\be$ norm
it is not true that for an arbitrary matrix function $\Phi\in L^\be(\mm_{m,n})$ there exists a dual extremal function in $H^1_0(\bS_1^{n,m})$. Moreover, it was shown in \cite{P3} that a dual extremal function exists if and only if the Hankel operator $H_\Phi:H^2(\C^n)\to H^2_-(\C^m)$ has a maximizing vector.

However, in the case $p=\be$, if a dual extremal function exists, then there exists a dual extremal function $\Psi$ such that $\rank\Psi(\z)=1$ almost everywhere on $\T$, see \cite{P3}.

\

\section{\bf Weird matrix functions}
\setcounter{equation}{0}

\

The main result of this section is a construction of a weird matrix function of size $2\times 2$.

\begin{lem}
\label{MC}
There exists a bounded $2\times2$ matrix function $B$ such that $B^*=B$, $\trace B(\z)=1$, $\z\in\T$,
the eigenvalues of
$B(\z)$, $\z\in\T$,  are positive and separated away from zero and there is no constant self-adjoint  matrix $C$ such that 
$$
\rank C=1\quad\mbox{and}\quad\trace B(\z)C=1, \quad\z\in\T.
$$
\end{lem}

\Pf Let $\a$ be a real bounded scalar functions, $\b$ a complex scalar bounded function such that the functions $\a$, $\b$, $\bar\b$,  and $\1$ 
are linearly independent, and the function $\a(1-\a)-|\b|^2$ is positive and separated away from zero. 
Put
$$
B=
\left(
\begin{array}{cc}
  \a&\b      \\[.2cm]
\bar\b  &1-\a       
\end{array}
\right).
$$
Clearly, $B^*=B$, the eigenvalues of $B(\z)$, $\z\in\T$, are positive and separated away from zero, and $\trace B(\z)=1$, $\z\in\T$. Suppose that $C$ is a self-adjoint constant matrix such that $\rank C=1$, and $\trace B(\z)C=1$, $\z\in\T$. Then $C$ has the form
$$
C=
\left(
\begin{array}{cc}
  a&b      \\[.2cm]
\bar b  & a^{-1}|b|^2  
\end{array}
\right),
$$
where $a$ is a nonzero real number and $b$ is a complex number. We have
$$
\trace B(\z)C=a\a(\z)+b\bar\b(\z)+\bar b\b(\z)+a^{-1}|b|^2(1-\a(\z))=1,\quad\z\in\T.
$$
Thus
$$
(a-a^{-1}|b|^2)\a(\z)+b\bar\b(\z)+\bar b\b(\z)+a^{-1}|b|^2-1=0.
$$
Since the functions  $\a$, $\b$, $\bar\b$, and $\1$ are linearly independent, this equality is impossible. $\bl$

Consider the Wiener--Masani factorization of $B$ (see \cite{WM}):
\bay
\label{mB}
B=\Psi^*\Psi,
\ey
where $\Psi$ is an invertible bounded analytic function in $\dd$. 
Put
$$
A=\Psi\Psi^*
$$
and consider the Wiener--Masani factorization of $A^2$:
$$
A^2=QQ^*.
$$
Let $U$ be the matrix function defined by
\bay
\label{dU}
U=\bar zQ^{-1}A.
\ey
Then $U$ is a unitary-valued function on $\T$:
$$
U^*U=A(Q^*)^{-1}Q^{-1}A=I.
$$
Clearly, 
$$
AU^{-1}=zQ\in H^\be_0(\mm_{2,2}).
$$
Let us show that $U$ is $p$-badly approximable. 

\begin{lem}
\label{U}
Let $A$ be a self-adjoint $2\times2$ matrix function such that
$$
\trace A(\z)=1,\quad\z\in\T,
$$
and the eigenvalues of $A(\z)$ are positive and separated away from zero.
Suppose that $U$ is a unitary-valued matrix function on $\T$ such that
$AU^{-1}\in H^\be_0(\mm_{2,2})$.
Then $U$ is a $p$-badly approximable matrix function.
\end{lem}



\Pf Let $F\in H^p(\mm_{2,2})$. For $\z\in\T$, we have
$$
\left|\trace\big((U-F)AU^*)(\z)\big)\right|\le\left\|(U-F)AU^*)(\z)\right\|_{\bS_1}
\le\|(U-F)(\z)\|_{\mm_{2,2}}.
$$
Thus by H\"older's inequality,
\begin{align*}
\|U-F\|_{L^p(\mm_{2,2})}&=\left(\int_\T \|U-F\|^p \,d\m\right)^{1/p}
\ge\left(\int_\T\left |\trace\big((U-F)AU^*\big)\right|^p \,d\m\right)^{1/p}\\[.2cm]
&\ge\int_\T \left|\trace\big((U-F)AU^*\big)\right| \,d\m
\ge\left|\int_\T \trace\big((U-F)AU^*\big)\,d\m\right|\\[.2cm]
&=\left|\int_\T\trace{(UAU^*)}\,d\m-\int_\T\trace{(FAU^*)}\,d\m\right|\\[.2cm]
&=\left|\int_\T\trace{(UAU^*)}\,d\m\right|=\left|\int_\T\trace A\,d\m\right|=1.
\end{align*}
Note that 
$$
\int_\T\trace{(FAU^*)}\,d\m=0,
$$
since $FAU^*\in H^p_0(\mm_{2,2})$. Thus $U$ is $p$-badly approximable. $\bl$


To prove that the matrix function $U$ defined by \rf{dU} is weird, we need the following lemma.

\begin{lem}
\label{Sp}
Let $A$ be a bounded positive definite matrix function on $\T$ whose inverse is also bounded
and let $A=\Psi\Psi^*$, where $\Psi$ is an invertible matrix function in $H^\be$.
A matrix function $F$ in $H^\be$ satisfies the equation
\begin{equation}
\label{uS}
AF^*=FA
\end{equation}
if and only if 
\begin{equation}
\label{par}
F=\Psi C\Psi^{-1}.
\end{equation}
where $C$ is a constant self-adjoint matrix.
\end{lem}

\Pf Put 
$$
C=\Psi^{-1}F\Psi.
$$
Then $C$ is an $H^\be$ matrix function and \rf{par} holds. By \rf{uS},  we have
$$
AF^*=A(\Psi^*)^{-1}C^*\Psi^*=FA=\Psi C\Psi^{-1}A.
$$
Since $A=\Psi\Psi^*$, we obtain
$$
\Psi\Psi^*(\Psi^*)^{-1}C^*\Psi^*=\Psi C\Psi^{-1}\Psi\Psi^*
$$
which implies $C=C^*$. Since $C$ is an $H^\be$ matrix function, it must be constant.

Clearly, if $C$ is a constant self-adjoint matrix and $F$ is defined by \rf{par}, then $F$ satisfies equation \rf{uS}.
$\bl$

\begin{thm}
\label{Wpb}
The matrix function $U$ defined by {\em\rf{dU}} is a weird $p$-badly approximable function.
\end{thm}

\Pf Assume that $U$ is respectable. By Lemma \ref{U}, $U$ is $p$-badly approximable.
Then $\|H_U\|=\|U\|_{L^p(\mm_{2,2})}=1$.

Let $f\in H^q(\C^2)$ be a maximizing vector of $H_U$ of norm 1. We have
$$
1=\|H_Uf\|_{L^2(\C^2)}=\|\pp_-Uf\|_{L^2(\C^2)}\le\|Uf\|_{L^2(\C^2)}=\|f\|_{L^2(\C^2)}\le\|f\|_{L^q(\C^2)}=1.
$$
Thus all inequalities in this chain of inequalities are equalities. 
The equality $\|f\|_{L^2(\C^2)}=\|f\|_{L^q(\C^2)}$ means that $\|f(\z)\|_{\C^2}=1$, $\z\in\T$, while the
equality
\lb$\|\pp_-Uf\|_{L^2(\C^2)}=\|Uf\|_{L^2(\C^2)}$ means that \mbox{$Uf\in H^2_-(\C^2)$}, and so
$Uff^*\in H^\be_-(\mm_{2,2})$ or, in other words, $ff^*U^{-1}\in H^\be_0(\mm_{2,2})$. 
Put
$$
F=ff^*A^{-1}.
$$
Then $F$ satisfies \rf{uS}. Hence, by Lemma \ref{Sp}, $F$ has the form $F=\Psi C\Psi^{-1}$, where $C$ is 
a constant self-adjoint matrix. Since $F$ has rank one on $\T$, it follows that $\rank C=1$. Clearly, 
$$
ff^*=FA=\Psi C\Psi^{-1}\Psi\Psi^*=\Psi C\Psi^*.
$$
Let $B=\Psi^*\Psi$ be the matrix function obtained in Lemma \ref{MC}. By \rf{mB}, we have
$$
\trace BC=\trace \Psi^*\Psi C=\trace \Psi C\Psi^*=\trace ff^*=1.
$$
This contradicts Lemma \ref{MC}. $\bl$

\medskip

{\bf Remark.} The results of Sections 4 and 5 show that the class $L^p(\mm_{m,n})$ splits into two subsets. The first subset consists of respectable matrix functions and for respectable matrix functions
$\Phi$ the distance $\dist_{L^p}\big(\Phi,H^p(\mm_{m,n})\big)$ can be computed by formula
\rf{fla}. The second subset  consists of weird matrix functions and for weird matrix functions $\Phi$
we have to find another formula to compute the distance $\dist_{L^p}\big(\Phi,H^p(\mm_{m,n})\big)$.
Such a formula will be obtained in the next section.

Let us explain that in a sense both the set of respectable matrix functions and the set of weird matrix functions are massive subsets of $L^p(\mm_{m,n})$. First of all, the set of weird matrix functions is open and nonempty, as we have just seen.

Secondly, if $\Phi\in L^p(\mm_{m,n})$ and $Q$ is an arbitrary function in $H^\be(\mm_{m,n})$, then
$\Phi$ is respectable if and only if $\Phi-Q$ is. Thus to characterize the set of respectable matrix functions, we can restrict ourselves to the case of $p$-badly approximable respectable matrix functions.
It is easy to see that the set of respectable badly approximable matrix functions has nonempty interior 
in the set of $p$-badly approximable matrix functions. Indeed, it is easy to verify that the $p$ badly approximable matrix function
\bay
\label{Phi}
\Phi=\left(\begin{matrix}\bar z&\0\\\0&\0\end{matrix}\right)
\ey
belongs to the interior of the set of respectable $p$-badly approximable functions.

However, we do not know whether the set of respectable matrix functions has nonempty interior in the space $L^p(\mm_{m,n})$. In particular, we do not know whether the matrix function $\Phi$ defined in \rf{Phi} belongs to the interior of the set of respectable matrix functions has nonempty interior in the space $L^p(\mm_{m,n})$.

\

\section{\bf Hankel operators on spaces of matrix-valued functions}
\setcounter{equation}{0}

\

We have already mentioned in the introduction that  the problem of analytic approximation of matrix functions can be reduced to the case of square matrix functions and beginning this section we assume that $\Phi\in L^p(\mm_{n})$.

For $\Phi\in L^p(\mm_{n})$, we consider the Hankel operator $\bs{H}_\Phi$
defined on the space $H^q(\bS_2^n)$ to the space $H^2_-(\bS_2^n)$ defined by
$$
\bs{H}_\Phi F=\pp_-\Phi F,
$$
where $\pp_-$ is an orthogonal projection from the space $L^2(\bS_2^n)$ onto the subspace
$H^2_-(\bS_2^n)\df L^2(\bS_2^n)\ominus H^2(\bS_2^n)$.

\begin{thm}
\label{HS}
Let \mbox{$\Phi\in L^p(\mm_{n})$}. Then
$$
\|\bs{H}_\Phi\|_{H^q(\bS_2^n)\to H^2_-(\bS_2^n)}=
\dist_{L^p}\big(\Phi, H^p(\mm_{n})\big).
$$
\end{thm}

\Pf Suppose that $\Phi\in L^p(\mm_{n})$. Then for $F\in H^q(\bS_2^n)$ and 
$Q\in H^p(\mm_{n})$, we have
\begin{align*}
\|\bs{H}_\Phi F\|_{L^2(\bS_2^n)}&=\|\pp_-\big((\Phi-Q)F\big)\|_{L^2(\bS_2^n)}
\le\|(\Phi-Q)F\|_{L^2(\bS_2^n)}\\[.2cm]
&\le\|\Phi-Q\|_{L^p(\mm_{n})}\|F\|_{H^q(\bS_2^n)}
\end{align*}
by H\"older's inequality.
Thus $\|\bs{H}_\Phi\|\le\dist_{L^p}\big(\Phi, H^p(\mm_{n})\big)$.

To prove the opposite inequality, we are going to use Theorem \ref{st} that has been deduced from
Sarason's theorem. Let $\Psi$ be a dual extremal function of $\Phi$, i.e., $\Psi$ belongs to 
$H^{p'}_0(\mm_{n})$ and satisfies \rf{de}. By Theorem \ref{st}, there exist matrix functions
\lb$F\in H^q(\mm_{n})$ and $G\in H^2_0(\mm_{n})$
such that
$$
\Psi=FG\quad\mbox{and}\quad
\|F\|_{L^q(\bS_2^n)}\|G\|_{L^2(\bS_2^n)}=1.
$$
Without loss of generality we may assume that $\|F\|_{L^q(\bS_2^n)}=1$ and $\|G\|_{L^2(\bS_2^n)}=1$.
We have
\begin{align*}
\|\bs{H}_\Phi\|&\ge|(\bs{H}_\Phi F,G^*)_{L^2(\bS_2^n)}|=
\left|\int_\T\trace\big((\bs{H}_\Phi F)G\big)\,d\m\right|\\[.2cm]
&=\left|\int_\T\trace(\Phi FG)\,d\m\right|
=\int_\T \trace\big(\Phi\Psi\big)\,d\m\\[.2cm]
&=\dist_{L^p}\big(\Phi,H^p(\mm_{m,n})\big)
\end{align*}
by \rf{de}. $\bl$

It follows immediately from Theorem \ref{HS} that 
$$
\|H_\Phi\|_{H^q(\C^n)\to H^2_-(\C^n)}\le\|\bs{H}_\Phi\|_{H^q(\bS_2^n)\to H^2_-(\bS_2^n)}.
$$
Note that this inequality can also be obtained easily from the definitions of $H_\Phi$ and $\bs{H}_\Phi$.

\begin{thm}
\label{un}
Let $\Phi\in L^p(\mm_{n})$ and let $Q$ be a best approximant to $\Phi$ in 
$H^p(\mm_{n})$. Then the following assertions hold:

{\em(i)} if $F\in H^q(\mm_{n})$ is a maximizing vector of $\bs{H}_\Phi$, then 
$(\Phi-Q)F\in H^2_-(\bS_2^n)$;

{\em(ii)}  the function
\bay
\label{dis}
\z\mapsto\|(\Phi-Q)(\z)\|_{\mm_{n}} 
\ey
does not depend on the choice of a best approximant $Q$; 

{\em(iii)} if $\bs{H}_\Phi$  has a maximizing vector $F$ such that $\rank F(\z)=n$ on a subset
of $\T$ of positive measure, then $\Phi$ has a unique best approximant in $H^p(\mm_{n})$;

{\em(iv)} if $F_1$ and $F_2$ are maximizing vectors of $\bs{H}_\Phi$, then 
$$
\|F_1(\z)\|_{\bS_2^n}=c\|F_2(\z))\|_{\bS_2^n}
$$
for some positive constant $c$;

{\em(v)} if $Q$ is a best approximant to $\Phi$ in 
$H^p(\mm_{n})$ and $F$ is a maximizing vector of $\bs{H}_\Phi$, then the matrix
$$
\frac1{\|(\Phi-Q)(\z)\|_{\mm_n}}(\Phi-Q)(\z),\quad\z\in\T,
$$
is isometric on the range of $F(\z)$.
\end{thm}

\Pf Let us fix a maximizing vector $F$ of $H_\Phi$. We have by H\"older's inequality,
\begin{align}
\label{ch}
\|\bs{H}_\Phi F\|_{L^2(\bS_2^n)}&=\|\pp_-\big((\Phi-Q)F\big)\|_{L^2(\bS_2^n)}
\le\|(\Phi-Q)F\|_{L^2(\bS_2^n)}\\[.2cm]
&\le\|\Phi-Q\|_{L^p(\mm_{n})}\|F\|_{H^q(\bS_2^n)}=
\|\bs{H}_\Phi\|_{H^q(\bS_2^n)\to H^2_-(\bS_2^n)}\|F\|_{H^q(\bS_2^n)}.\nonumber
\end{align}
Since $\|\bs{H}_\Phi F\|_{L^2(\bS_2^n)}=\|\bs{H}_\Phi\|\cdot\|F\|_{H^q(\bS_2^n)}$,
it follows that both inequalities in \rf{ch} are equalities. 

The fact that the first inequality in \rf{ch} turns into equality means that
\lb$(\Phi-Q)F\in H^2_-(\bS_2^n)$, i.e.,
\begin{equation}
\label{Q}
(\Phi-Q)F=\bs{H}_\Phi F
\end{equation}
which proves (i).

To prove (iii), we observe that since $F\in H^q(\bS_2^n)$, it follows that if $\rank F(\z)=n$ on a 
set of positive measure, then
$\rank F(\z)=n$,  $\z\in\T$, almost everywhere on $\T$. Hence,
$$
\Phi-Q=(\bs{H}_\Phi F)F^{-1},
$$
and so $Q$ is uniquely determined by $\Phi$.

The fact that the second inequality in \rf{ch} turns into equality means that there exists $c>0$ such that
\bay
\label{uni}
\|(\Phi-Q)(\z)\|_{\mm_{n}}^p=c\|F(\z)\|_{\bS_2^n}^q,\quad\z\in\T,
\ey
and
\bay
\label{mv}
\|(\Phi-Q)(\z)F(\z)\|_{\bS_2^n}=\|(\Phi-Q)(\z)\|_{\mm_{n}}\|F(\z)\|_{\bS_2^n},\quad\z\in\T.
\ey
Clearly,  (iv) follows immediately from \rf{uni}. 

If we normalize the maximizing vector $F$ by the condition
\bay
\label{norm}
\|F\|_{L^q(\bS_2)}^q=\|\bs{H}_\Phi\|^p,
\ey
then integrating \rf{uni}, we obtain
$$
\|\bs{H}_\Phi\|^p=\|\Phi-Q\|_{L^p(\mm_{n})}^p=c\|F\|_{L^q(\bS_2)}^q.
$$
Hence, under condition \rf{norm}, 
\bay
\label{normal}
\|(\Phi-Q)(\z)\|_{\mm_{n}}^p=\|F(\z)\|_{\bS_2^n}^q,\quad\z\in\T,
\ey
and so the function \rf{dis} is uniquely determined by $\Phi$. This proves (ii).

It remains to observe that (v) follows from \rf{mv} and from the fact that for $n\times n$ matrices
$A$ and $B$ the equality
$$
\|AB\|_{\bS_2^n}=\|A\|\cdot\|B\|_{\bS_2^n}
$$
holds if and only if the restriction of $A$ to the range of $B$ is a multiple of an isometry.
$\bl$

\medskip

{\bf Definition.} For a function $\Phi\in L^p(\mm_{n})$, the function \rf{dis} is called
the {\it distance function} of $\Phi$. We denote the distance function of $\Phi$ by
${\rm d}_\Phi$:
\bay
\label{df}
{\rm d}_\Phi(\z)=\|(\Phi-Q)(\z)\|_{\mm_{n}},\quad\z\in\T,
\ey
where $Q$  is an arbitrary best approximant to $\Phi$.

\medskip

The following result describes the set of all nonzero distance functions of 
matrix functions in $L^p(\mm_{n})$.

\begin{thm}
\label{log}
Let $d\ge0$ be a nonzero function in $L^p$. Then $d$ is
the distance function of a matrix function $\Phi\in L^p(\mm_{n})$
if and only if $\log d\in L^1$.
\end{thm}

\Pf If $\Phi\in L^p(\mm_n)\setminus H^p(\mm_n)$ and $Q$ is a $p$-best approximant to $\Phi$
and $d(\z)=\|(\Phi-Q)(\z)\|_{\mm_n}$, the fact that 
$\log d\in L^1$ follows immediately from \rf{uni}.
 
The converse follows from Corollary \ref{mod}
by considering matrix functions of the form
$$
\left(\begin{matrix}\f&\0&\cdots&\0\\
\0&\0&\cdots&\0\\
\vdots&\vdots&\ddots&\vdots\\
\0&\0&\cdots&\0
\end{matrix}
\right).
\quad\bl
$$

\begin{thm}
\label{rcr}
$\Phi\in L^p(\mm_{n})\setminus H^p(\mm_{n})$. Then $\Phi$ is respectable if and only if there exists a maximizing vector 
$F$ of $\bs{H}_\Phi$ such that 
\bay
\label{r1}
\rank F(\z)=1,\quad\z\in\T.
\ey
\end{thm}

\Pf Suppose that $\Phi$ is respectable. Consider the Hankel operator
$$
H_\Phi:H^q(\C^n)\to H^2_-(\C^n).
$$
Let $f\in H^q(\C^n)$ be a maximizing vector of $H_\Phi$. Define the matrix function
\lb$F\in H^q(\mm_{n})$ by
$$
F=\left(\begin{array}{cccc}f&\0&\cdots&\0\end{array}\right).
$$
It is obvious that $\rank F(\z)=1$ for $\z\in\T$. Clearly, 
$$
\|F(\z)\|_{\bS_2}=\|f(\z)\|_{\C^n}\quad\mbox{and}\quad 
\|(\bs{H}_\Phi F)(\z)\|_{\bS_2^n}=\|(H_\Phi f)(\z)\|_{\C^n},
$$
and so $F$ is a maximizing vector of $\bs{H}_\Phi$. 

To prove the converse, we may assume that $\|H_\Phi\|=1$.
Suppose that $F$ is a maximizing vector of $\bs{H}_\Phi$ of norm 1 that satisfies \rf{r1}.
Let $Q$ be a best approximant to $\Phi$ in $H^p(\mm_{n})$. By Theorem \ref{un}, we have
$(\Phi-Q)F\in H^2_-(\mm_{n})$. Put
$$
G=\frac{1}{\|\bs{H}_\Phi\|}\Big((\Phi-Q)F\Big)^*\quad\mbox{and}\quad\Psi=FG\in H^{p'}_0(\mm_{n}).
$$
Clearly, 
$$
\rank\Psi(\z)=1,\quad\z\in\T,\quad\mbox{and}\quad\|\Psi\|_{L^{p'}(\bS_1^n)}\le1.
$$
Let us show that $\Psi$ is a dual extremal function of $\Phi$. Assuming that \rf{norm} holds, we have by
\rf{normal},
\begin{align*}
\int_\T\trace(\Phi\Psi)\,d\m&=\int_\T\trace\big((\Phi-Q) FG\big)\,d\m\\[.2cm]
&=\int_\T\trace\Big(\pp_-\big((\Phi-Q) F\big)G\Big)\,d\m\\[.2cm]
&=\int_\T\trace\big((\bs{H}_\Phi F)G\big)\,d\m
=(\bs{H}_\Phi F,G^*)_{L^2(\bS_2^n)}\\[.2cm]
&=\frac{1}{\|\bs{H}_\Phi\|}\big\|\bs{H}_\Phi F\big\|^2_{L^2(\bS_2^n)}=\|\bs{H}_\Phi\|=1.
\end{align*}
Thus $\Psi$ is a dual extremal function of rank 1, and so by Theorem \ref{Rf},
$\Phi$ is respectable. $\bl$

Note that the computation, in fact, shows that $\|\Psi\|_{L^{p'}(\bS_1^n)}=1$.

\begin{cor}
\label{edi}
Let $\Phi$ be a weird function in $L^p(\mm_{2})\setminus H^p(\mm_{2})$. 
Then $\Phi$ has a unique best approximant in $H^p(\mm_{2})$.
\end{cor}

\Pf By Theorem \ref{un}, if $\Phi$ has more than one best approximant, then each maximizing vector 
$F$ of $\bs{H}_\Phi$ has rank 1 almost everywhere on $\T$. By Theorem \ref{rcr}, the function
$\Phi$ is respectable. $\bl$

We consider now for a function $\Phi\in L^p(\mm_{n})$,  the family of Hankel operators
$H_\Phi^{\{k\}}: H^q\big(\bS_2^{n,k}\big)\to H^2_-\big(\bS_2^{n,k}\big)$, $1\le k\le n$, defined by
$$
H_\Phi^{\{k\}}F=\pp_-\Phi F,\quad F\in  H^q(\bS_2^{n,k}).
$$
Clearly, $H_\Phi^{\{1\}}=H_\Phi$ and $H_\Phi^{\{n\}}=\bs{H}_\Phi$.

\begin{thm}
\label{rankk}
Let $\Phi\in L^p(\mm_{n})\setminus H^p(\mm_{n})$ and let $1\le k\le n$. The following are equivalent:

{\em(i)} there exists a maximizing vector $F$ of $\bs{H_\Phi}$ such that
\bay
\label{mvk}
\rank F(\z)\le k,\quad\z\in\T;
\ey

{\em(ii)} the following distance formula holds:
$$
\left\|H_\Phi^{\{k\}}\right\|=\dist_{L^p}\big(\Phi,H^p(\mm_{n})\big).
$$
\end{thm}

Note that a standard argument with analyticity properties of minors shows that $\rank F(\z)$ is constant 
for almost all $\z$ in $\T$.

\medskip

{\bf Proof of Theorem \ref{rankk}.} Suppose that (ii) holds. Let $G\in H^q(\bS_2^{n,k})$ be a maximizing vector of 
${H}^{(k)}_\Phi$ (observe that ${H}^{(k)}$ is compact). Consider the matrix function $F\in H^q(\bS_2^n)$ obtained from $G$ by adding $n-k$ zero columns. Clearly.
$$
\|\bs{H}_\Phi F\|_{L^2(\bS_2^n)}=\left\|H_\Phi^{\{k\}}G\right\|_{L^2(\bS_2^{n,k})}
=\dist_{L^p}\big(\Phi,H^p(\mm_{n})\big)\|F\|_{L^q(\bS_2^n)}.
$$ 
Thus $F$ is a maximizing vector of $\bs{H}_\Phi$ that satisfies \rf{mvk}.

Suppose now that $F$ is a maximizing vector of $\bs{H}_\Phi$ such that
$$
\rank F(\z)=k,\quad\z\in\T.
$$
Without loss of generality we may assume that $\|F\|_{L^q(\bS_2^n)}=1$.
As in the proof of Theorem \ref{rcr}, consider a best approximant $Q$ to $\Phi$ in $H^p(\mm_{n})$
and define the matrix functions $G$ and $\Psi$ as in that proof.
Then
$$
\rank\Psi(\z)=k,\quad\z\in\T,\quad\mbox{and}\quad\|\Psi\|_{L^{p'}(\mm_{n})}=1.
$$
The fact that $\Psi$ is a dual extremal function of $\Phi$ can be verified as in the proof of Theorem 
\ref{rcr}.

By Theorem \ref{nkq2}, $\Psi$ admits a factorization $\Psi=F_\flat G_\flat$, where
$F_\flat \in H^q(\mm_{n,k})$, $G_\flat\in H^2_0(\mm_{k,n})$ and
$\|\Psi\|_{L^{p'}(\bS_1^n)}=\|F_\flat\|_{L^q(\bS_2^{n,k})}\|G_\flat\|_{L^2(\bS_2^{k,n})}$.

We claim that $\|\bs{H}_\Phi\|=\left\|H_\Phi^{\{k\}}\right\|$ and $F_\flat$ is a maximizing vector of
$H_\Phi^{\{k\}}$. This can be proved in  the same way as in the proof of Theorem 
\ref{HS}. Indeed,
without loss of generality we may assume that $\|F_\flat\|_{L^q(\bS_2^{n,k})}=1$ and 
$\|G_\flat\|_{L^2(\bS_2^{k,n})}=1$.
Then
\begin{align*}
\|\bs{H}_\Phi\|&\ge
\left\|H^{\{k\}}_\Phi\right\|\ge\left|\left(H^{\{k\}}_\Phi F_\flat,G_\flat^*\right)_{L^2(\bS_2^n)}\right|=
\left|\int_\T\trace\Big(\big(H^{\{k\}}_\Phi F_\flat\big)G_\flat\Big)\,d\m\right|\\[.2cm]
&=\left|\int_\T\trace\big(\pp_-(\Phi F_\flat),G_\flat\big)\,d\m\right|
=\left|\int_\T\trace(\Phi F_\flat G_\flat)\,d\m\right|\\[.2cm]
&=\int_\T \trace\big(\Phi\Psi\big)\,d\m
=\dist_{L^p}\big(\Phi,H^p(\mm_{n})\big)=\|\bs{H}_\Phi\|
\end{align*}
by \rf{de} and Theorem \ref{HS}. $\bl$

\medskip 

{\bf Definition.} A matrix function $\Phi\in L^p(\mm_{n})$ is said to {\it have order $k$} if $k$ is the smallest number such that 
$$
\dist_{L^p}\big(\Phi,H^p(\mm_{n})\big)=\left\|H_\Phi^{\{k\}}\right\|.
$$

\medskip

Clearly, a matrix function $\Phi$ is respectable if and only if it has order 1.

The reasoning given in the proof of Theorem \ref{rankk} allows us to obtain the following formulae for the
order of a matrix function in $L^p$.

\begin{thm}
\label{reas}
Let $\Phi$ be a matrix function in $L^p(\mm_n)$. Then the following assertion hold:

{\em(i)} the order of $\Phi$ is the minimal number $k$, for which there exists a maximizing vector $F$ of
$\bs{H}_\Phi$ that satisfies {\em\rf{mvk}}.

{\em(ii)} the order of $\Phi$ is the minimal number $k$ such that $\Phi$ has a dual extremal function $\Psi$
satisfying
$$
\rank\Psi(\z)\le k,\quad\z\in\T.
$$
\end{thm}

\Pf It is easy to see that the proof of Theorem \ref{reas} is contained in the proof of Theorem \ref{rankk}.
$\bl$

In \S 7 we obtain one more formula for the order of $\Phi$, see Theorem \ref{ord}.

We can obtain now an analog of Theorem \ref{un} for the Hankel operators $H_\Phi^{\{k\}}$.

\begin{thm}
\label{anal}
Let $\Phi$ be a matrix function in $L^p(\mm_n)$ such that
$$
\dist_{L^p}\big(\Phi,H^p(\mm_n)\big)=\left\|H_\Phi^{\{k\}}\right\|.
$$
Then the following assertions hold:

{\em(i)} if $F\in H^q(\mm_{n,k})$ is a maximizing vector of $H_\Phi^{\{k\}}$, then 
$(\Phi-Q)F\in H^2_-(\bS_2^{n,k})$;

{\em(ii)} if $F_1$ and $F_2$ are maximizing vectors of $H^{(k)}_\Phi$, then 
$$
\|F_1(\z)\|_{\bS_2^n}=c\|F_2(\z))\|_{\bS_2^n}
$$
for some positive constant $c$;

{\em(iii)} if $Q$ is a best approximant to $\Phi$ in 
$H^p(\mm_{n})$ and $F$ is a maximizing vector of $H_\Phi^{\{k\}}$, then the matrix
$$
\frac1{\|(\Phi-Q)(\z)\|}(\Phi-Q)(\z),\quad\z\in\T,
$$
is isometric on the range of $F(\z)$.
\end{thm}

Theorem \ref{anal} can be proved in the same way as Theorem \ref{un}.

\medskip

{\bf Remark.} Note that in the case $p=\be$ and $k=1$, (ii) is very far from being true. Indeed, we can take two different scalar outer functions $h_1$ and $h_2$ in $H^2$ and consider the matrix function
$\Phi$ defined by
$$
\Phi=
\left(\begin{matrix}\bar z\frac{\bar h_1}{h_1}&\0\\[.2cm]\0&\bar z\frac{\bar h_2}{h_2}\end{matrix}\right).
$$
It is easy to see that $\Phi$ is badly approximable, $\|H_\Phi\|_{H^2(\C^2)\to H^2_-(\C^2)}=1$, and the vector functions 
$$
f_1=\left(\begin{matrix}h_1\\\0\end{matrix}\right)\quad\mbox{and}\quad
f_2=\left(\begin{matrix}\0\\h_1\end{matrix}\right)
$$
are maximizing vectors of $H_\Phi$, though the functions
$$
\z\mapsto\|f_1(\z)\|_{\C^2}=|h_1(\z)|\quad\mbox{and}\quad
\z\mapsto\|f_2(\z)\|_{\C^2}=|h_2(\z)|
$$
do not have to be proportional.

\

\section{\bf  $\bs{p}$-badly approximable functions}
\setcounter{equation}{0}

\

In this section we characterize the set of all badly approximable functions in terms of certain special factorizations. Such factorizations allow us in this section to obtain a parametrization of all $p$-best approximants to a given matrix function in $L^p(\mm_{n})$ in the case when such a best approximant is not unique.

To describe the set of $p$-badly approximable matrix functions, we prove the following result
that can be considered as an analog of the corresponding result for analytic approximation in the 
$L^\be$ norm, see \cite{P}, Ch. 14, \S\,15.

\begin{thm}
\label{ptf}
Let $\Phi\in L^p(\mm_{n})$ and let $Q$ be a best approximant to $\Phi$ in 
$H^p(\mm_{n})$. Then $\Phi-Q$ admits the following factorization
\bay
\label{ff}
\Phi-Q=\W^*\left(\begin{array}{cc}\D&\0\\\0&\Phi_\#\end{array}\right)\V^*, 
\ey
where $\V$ and $\W^{\rm t}$ are $k$-balanced matrix functions for some $k\le n$,
$\D$ is a  \lb$k\times k$ $p$-badly approximable matrix function such that the matrix function 
${\rm d}_\Phi^{-1}\D$ is unitary-valued, and $\Phi_\#$ is a matrix function such that 
$$
\|\Phi_\#(\z)\|_{\mm_{n-k}}\le\|\D(\z)\|_{\mm_{k}}, \quad\z\in\T.
$$ 
\end{thm}

\Pf Clearly, without loss of generality we may assume that $Q=\0$, i.e., $\Phi$ is a $p$-badly approximable matrix function.

Suppose that $\bs{H}_\Phi$ has a maximizing vector of rank $k$. In the proof of Theorem 
\ref{rankk} we have shown that $\left\|H^{\{k\}}_\Phi\right\|=\dist_{L^p}\big(\Phi,H^p(\mm_{n})\big)$
and there exists a maximizing vector $F\in H^q(\mm_{n,k})$ of $H_\Phi^{\{k\}}$ such that 
$\rank F(\z)=k,\quad\z\in\T$. Consider the inner-outer factorization of $F^{\rm t}$:
$$
F^{\rm t}=\cO_1^{\rm t}F_{\rm co}^{\rm t}.
$$
Then
$$
F=F_{\rm co}\cO_1
$$
It is easy to see that $\cO_1$ is an inner matrix function of size $k\times k$ and $F_{\rm co}$
is a co-outer matrix function of size $n\times k$. It follows easily from (i) of Theorem  \ref{anal} 
that $F_{\rm co}$ is a maximizing vector 
of $H_\Phi^{\{k\}}$. Without loss of generality we may thus assume that $F$ is co-outer.

Let $G$ be the function in $H^2(\mm_{n,k})$ defined by
$$
G(\z)=\ov{\z}\,\ov{\big(H_\Phi^{\{k\}}F\big)(\z)}.
$$
By (i) and (iii) of Theorem \ref{anal}, we know that $G=\bar z\ov{\Phi F}$ has rank $k$ on $\T$.

Similarly, we can consider the inner-outer factorization of $G^{\rm t}$ and obtain a factorization
$$
G=G_{\rm co}\cO_2,
$$
where $\cO_2$ is an inner matrix function of size $k\times k$ and $G_{\rm co}$ is a co-outer
matrix function.

Consider now the inner-outer factorization of $F$
$$
F=\U F_{\rm o}.
$$
Since $\rank F(\z)=k$ almost everywhere on $\T$, it is easy to see that $\U$ has size $n\times k$.
Similarly, we can consider the inner-outer factorization of $G_{\rm co}$:
$$
G_{\rm co}=\O G_{\rm o}
$$
and $\O$ has size $n\times k$.

We can consider now balanced completions $\V$ and $\W$ of $\U$ and $\O$:
\bay
\label{cc}
\V=\left(\begin{array}{cc}\U&\ov{\Theta}\end{array}\right)\quad\mbox{and}\quad
\W^{\rm t}=\left(\begin{array}{cc}\O&\ov{\Xi}\end{array}\right),
\ey
where $\Theta$ and $\Xi$ are inner and co-outer matrix functions such that the matrix functions
$\V$ and $\W$  defined by \rf{cc} are unitary-valued (see \S\,2).

Let 
$$
A=\W\Phi\V.
$$
By Theorem \ref{anal}, $H^{\{k\}}_\Phi F=\Phi F$, and so
$$
\Phi F=\W^*A\V^*F=
\W^*A\left(\begin{array}{c}\U^*\\\Theta^{\rm t}\end{array}\right)\U F_{\rm o}
=\W^*A\left(\begin{array}{c}F_{\rm o}\\\0\end{array}\right)=\bar z\ov{G}=\bar z\ov{G_{\rm co}\cO_2}.
$$
Thus
$$
A\left(\begin{array}{c}F_{\rm o}\\\0\end{array}\right)=
\W\Phi F=\bar z\left(\begin{array}{c}\O^{\rm t}\\\Xi^*\end{array}\right)\ov{\O G_{\rm   o}\cO_2}
=\left(\begin{array}{c}\bar z \ov{G_{\rm o}\cO_2}\\\0\end{array}\right).
$$
Clearly, $\|A(\z)\|_{\mm_{n}}=\|\Phi(\z)\|_{\mm_{n}}$, $\z\in\T$, and by Theorem \ref{un} (see \rf{mv}),
$\left(\begin{array}{c}F_{\rm o}(\z)\\0\end{array}\right)$ is a maximizing vector of $A(\z)$ for almost all $\z\in\T$. Let 
$$
A=\left(\begin{array}{cc}A_{11}&A_{12}\\A_{21}&A_{22}\end{array}\right),
$$
where $A_{11}$ has size $k\times k$.

By Theorem \ref{anal}, the matrices $\|A(\z)\|^{-1}A_{11}(\z)$ take unitary values almost everywhere on $\T$.
It is easy to verify (see e.g., \cite{P}, Lemma 15.5 of Ch.~14) that $A_{21}=\0$, $A_{12}=\0$,
and $\|A_{22}(\z)\|\le\|A_{11}(\z)\|$, $\z\in\T$. 

Clearly, $\|A(\z)\|={\rm d}_\Phi(\z)$. Put $\D=A_{11}$. Then
$({\rm d}_\Phi)^{-1}\D$ is a unitary-valued matrix function and
$$
\Phi=\W^*\left(\begin{array}{cc}\D&\0\\\0&\Phi_\#\end{array}\right)\V^*,
$$
where $\Phi_\#\df A_{22}$. Obviously, $\|\Phi_\#(\z)\|\le\rm{d_\Phi}(\z)$,
$\z\in\T$.

It is easy to see that  
$\bs{H}_\D F_{\rm o}=\bar z\ov{G_{\rm o}\cO_2}$, and so
$\|\bs{H}_\D\|=\|\bs{H}_\Phi\|$, which implies that 
$\D$ is a $p$-badly approximable matrix function. $\bl$

\medskip

{\bf Remark 1.} Note that the matrix function $\D$ is determined by the choice of a maximizing vector
and it does not depend on the choice of a $p$-best approximant $Q$. It is also clear that the $k$-balanced matrix functions $\V$
and $\W$ do not depend on the choice of $Q$ either.

\medskip

{\bf Remark 2.} Clearly, we can always take $k$ to be the order of $\Phi$. However, the choice of $k$ is not always unique. For example, if $\f\in L^p$ is a scalar $p$-badly approximable function and 
$f\in H^q$ is a maximizing vector of $H_\f$, then it is easy to see that 
$\Phi=\left(\begin{array}{cc}\f&\0\\\0&\f\end{array}\right)$ is a respectable $p$-badly approximable 
matrix function  and the matrix function
$F=\left(\begin{array}{cc}f&\0\\\0&f\end{array}\right)$ is a maximizing vector
of $\bs{H}_\Phi=H_\Phi^{\{2\}}$. Thus the matrix function $\Phi$ admits factorizations of the form 
\ref{ptf} with $k=1$ and $k=2$.

\medskip

{\bf Definition.}
We say that a matrix function $\Phi\in L^p(\mm_{n})$ {\it has gender} $k$ if $k$ is the maximal number such that $\bs{H}_\Phi$ has a maximizing vector of rank $k$. Clearly, in Theorem \ref{ptf} we can take $k$ to be the gender of $\Phi$.

\medskip

Factorizations of the form \rf{ff} allow us to obtain one more formula for the order of matrix functions in $L^p$.

\begin{thm}
\label{ord}
Let $\Phi$ be a matrix function in $L^p(\mm_n)$ and let $Q$ is a $p$-best approximant to $\Phi$.
Then  the order of $\Phi$ is the minimal number $k$ such that $\Phi-Q$ admits a factorization as in
{\em \rf{ff}} with $k$-balanced matrix functions $\V$ and $\W^{\rm t}$.
\end{thm}

\Pf The proof of Theorem \ref{ptf} shows that 
if $k$ is the order of $\Phi$, then $\Phi-Q$ admits a factorization of the form \rf{ff} with
$k$-balanced matrix functions $\V$ and $\W^{\rm t}$.

Suppose now that \rf{ff} holds with
$k$-balanced matrix functions $\V$ and $\W^{\rm t}$. Suppose that $\V$ and $\W$ are given by \rf{cc}.

Let $G\in H^q(\bS_2^k)$ be a maximizing vector of $\bs{H}_\D$. Consider the matrix function $F\in H^q(\bS_2^{n,k})$
defined by
$$
F=\U G.
$$
We have
\begin{align*}
(\Phi-Q)F&=\W^*\left(\begin{matrix}\D&\0\\\0&\Phi_\#\end{matrix}\right)
\left(\begin{matrix}\U^*\\\Theta^{\rm t}\end{matrix}\right)\U G\\[.2cm]
&=\W^*\left(\begin{matrix}\D&\0\\\0&\Phi_\#\end{matrix}\right)
\left(\begin{matrix}G\\\0\end{matrix}\right)
=\left(\begin{matrix}\ov{\O}&\Xi\end{matrix}\right)
\left(\begin{matrix}\D G\\\0\end{matrix}\right).
\end{align*}
Since $G$ is a maximizing vector of $\bs{H}_\D$ and $\D$ is a $p$-badly approximable matrix
function, it follows from Theorem \ref{un} that $\D G=\bs{H}_\D G$, and so
\begin{align*}
(\Phi-Q)F=\ov{\O}\bs{H}_\D G.
\end{align*}
It is easy to see that $F$ is a maximizng vector of $H_\Phi^{\{k\}}$ and $\|\bs{H}_\Phi\|=\big\|H_\Phi^{\{k\}}\big\|$.
This proves the result. $\bl$

\medskip

{\bf Remark.} If $\Phi-Q$ admits a factorization as in \rf{ff} with $k$-balanced $\V$ and $\W$, where $k$ is the order of $\Phi$,
and $\L$ is a dual extremal function of $\D$, then a dual extremal function $\Psi$ of $\Phi$ satisfying the condition
$$
\rank\Psi(\z)=k,\quad\z\in\T,
$$
can be obtained by the following explicit formula:
$$
\Psi=\U\L\O^{\rm t}.
$$

Indeed, we may assume without loss of generality that $\Phi$ is $p$-badly approximable and
$Q=\0$. We have
$$
\|\Psi\|_{L^{p'}(\bS_1^n)}=\|\L\|_{L^{p'}(\bS_1^k)}=1.
$$
Clearly,
\begin{align*}
\trace\big(\Phi(\z)\Psi(\z)\big))&=
\trace\left(\O^{\rm t}(\z)
\Phi(\z)
\U(\z)\L(\z)
\right)\\[.2cm]
&=\trace\left(
\left(\begin{matrix}I_k&0\end{matrix}\right)
\left(\begin{matrix}\D(\z)&0\\0&\Phi_\#(\z)\end{matrix}\right)
\left(\begin{matrix}I_k\\0\end{matrix}\right)
\L(\z)\right)\\[.2cm]
&=\trace\Big(\D(\z)\L(\z)\Big).
\end{align*}
Thus
\begin{align*}
\int_\T\trace(\Phi(\z)\Psi(\z))\,d\m(\z)&=\int_\T\trace\big(\D(\z)\L(\z)\big)\,d\m(\z)\\[.2cm]
&=\|\D\|_{L^p(\mm_k)}=\|\Phi\|_{L^p(\mm_n)}.\quad\bl
\end{align*}

\medskip

Let us now describe all $p$-badly approximable matrix functions. Note that similar results hold in the case $p=\be$
under certain restrictions on the function, see \cite{P}, Ch. 14, \S\,15.

\begin{thm}
\label{dba}
Let $\Phi$ be matrix function in 
$L^p(\mm_{n})$. Then $\Phi$ is $p$-badly approximable if and only if there exists $k\le n$ such that
$\Phi$ admits a factorization
\bay
\label{pbaf}
\Phi=\W^*\left(\begin{array}{cc}\D&\0\\\0&\Phi_\#\end{array}\right)\V^*, 
\ey
where $\V$ and $\W^{\rm t}$ are $k$-balanced matrix functions,
$\D$ is a  $k\times k$ $p$-badly approximable matrix function such that the matrix function 
${\rm d}_\Phi^{-1}\D$ is unitary-valued, and $\Phi_\#$ is a matrix function such that 
$\|\Phi_\#(\z)\|_{\mm_{n-k}}\le\|\D(\z)\|_{\mm_{k}}$
for almost all $\z\in\T$.
\end{thm}

\Pf The fact that $p$-badly approximable matrix functions $\Phi$ admit factorizations of the form 
\rf{pbaf} follows immediately from Theorem \ref{ptf}.

Suppose now that $\Phi$ is given by \rf{pbaf}. Consider the Hankel operator
\lb$\bs{H}_\D:H^q(\bS_2^k)\to H^2_-(\bS_2^k)$. Let $\F\in H^q(\bS_2^k)$
be a maximizing vector of $\bs{H}_\D$.
Since $\D$
is $p$ badly approximable, it follows from Theorems \ref{HS} and \ref{un} that 
$$
\bs{H}_\D\F=\D\F\quad \mbox{and}\quad
\|\bs{H}_\D\F\|_{L^2(\bS_2^k)}=\|\Phi\|_{L^p(\mm_{n})}\|\F\|_{L^q(\bS_2^k)}.
$$
Consider the matrix function $F=\U\F$, where $\U$ is as in \rf{cc}. We have
\begin{align*}
\Phi F&=\W^*\left(\begin{array}{cc}\D&\0\\\0&\Phi_\#\end{array}\right)
\left(\begin{array}{c}\U^*\\\Theta^{\rm t}\end{array}\right)\U\F\\[.2cm]
&=\W^*\left(\begin{array}{cc}\D&\0\\\0&\Phi_\#\end{array}\right)
\left(\begin{array}{c}\F\\\0\end{array}\right)\\[.2cm]
&=\left(\begin{array}{cc}\ov{\O}&\Xi\end{array}\right)
\left(\begin{array}{c}\D\F\\\0\end{array}\right)=
\ov{\O}\D\F\in H^2_-(\bS_2^{n,k}).
\end{align*}
Thus
\begin{align*}
\left\|H_\Phi^{\{k\}} F\right\|_{L^2(\bS_2^{n,k})}&=\|\D\F\|_{L^2(\bS_2^k)}
=\|\bs{H}_\D\F\|_{L^2(\bS_2^k)}\\[.2cm]
&=\|\Phi\|_{L^p(\mm_{n})}\|\F\|_{L^2(\bS_2^k)}=\|\Phi\|_{L^p(\mm_{n})}\|F\|_{L^2(\bS_2^{n,k})}.
\end{align*}
It follows that $\left\|H_\Phi^{\{k\}}\right\|=\|\Phi\|_{L^p(\mm_{n})}$, and so $\Phi$ is
$p$-badly approximable. $\bl$

The next result allows us to parametrize all best approximants in the case when there are
more than one best approximant. A similar result also holds in the case of approximation in the norm of $L^\be$
under certain restrictions on $\Phi$, see \cite{P}, Ch. 14, \S\,15.

\begin{thm}
\label{paf}
Let $\Phi$ and $Q$ be as in Theorem {\em\ref{ptf}} and let $\Phi-Q$ be factorizred as
in {\em\rf{ff}}. A matrix function $R\in H^p(\mm_{n})$ is a $p$-best approximant to $\Phi$ if and only if there exists a matrix function $R_\#\in H^p(\mm_{n-k})$ such that
\bay
\label{ppba}
\Phi-R=\W^*\left(\begin{array}{cc}\D&\0\\\0&\Phi_\#-R_\#\end{array}\right)\V^*, 
\ey
and
\bay
\label{nerv}
\|(\Phi_\#-R_\#)(\z)\|_{\mm_{n-k}}\le\|\D(\z)\|_{\mm_{k}}, \quad\z\in\T.
\ey
\end{thm}

We need the following lemma.

\begin{lem}
\label{VW}
Let $\V$ and $\W^{\rm t}$ be $k$-balanced matrix functions of size $n\times n$. Then
$$
\W H^p(\mm_{n})\V\bigcap
\left(\begin{array}{cc}\0&\0\\\0&L^p(\mm_{n-k})
\end{array}\right)
=\left(\begin{array}{cc}\0&\0\\\0&H^p(\mm_{n-k})
\end{array}\right).
$$
\end{lem}

For $p=\be$ this is Theorem 1.8 of Ch. 14 of \cite{P}. The proof given in \cite{P} also works in our case.

\medskip

{\bf Proof of Theorem \ref{paf}.} Suppose that $R$ is a best approximant to $\Phi$. Then by Theorem \ref{ptf}, $\Phi-R$ admits
a factorization 
$$
\Phi-R=\W^*\left(\begin{array}{cc}\D&\0\\\0&\Phi_\flat\end{array}\right)\V^*,
$$
where $\D$, $\V$, and $\W$ are as in \rf{ptf} and $\Phi_\flat$ is a matrix function such that 
\lb$\|\Phi_\flat(\z)\|_{\mm_{n-k}}\le\|\D(\z)\|_{\mm_{k}}$, $\z\in\T$. Then
$$
R-Q=\W^*\left(\begin{array}{cc}\0&\0\\\0&\Phi_\#-\Phi_\flat\end{array}\right)\V^*.
$$
By Lemma \ref{VW}, $R_\#\df\Phi_\#-\Phi_\flat\in H^p(\mm_{n-k})$. 

Conversely, suppose that $R_\#$ is a matrix function in $H^p(\mm_{n-k})$ such  that 
\rf{nerv} holds.
Then by Lemma \ref{VW}, there exists 
$R\in H^p(\mm_{n})$ such that
$$
R-Q=\W^*\left(\begin{array}{cc}\0&\0\\\0&R_\#\end{array}\right)\V^*.
$$
Then \rf{ppba} holds. It follows easily from Theorem \ref{dba} that $R$ is a $p$-best analytic approximant to $\Phi$.
$\bl$ 

Theorem \ref{paf} says that to describe all $p$ best approximants, we should describe all functions
$R_\#\in H^p(\mm_{n-k})$ such that \rf{nerv} holds.
By Theorem \ref{log} there exists a scalar outer function $\vk$ in $H^p$ such that 
$|\vk(\z)|=\|\D(\z)\|_{\mm_{k}}={\rm d}_\Phi(\z)$, $\z\in\T$. Clearly, a matrix function $R_\#$
in $H^p(\mm_{n-k})$ satisfies \rf{nerv} if and only if the matrix function $\vk^{-1}R_\#$
satisfies the inequality
$$
\|\vk^{-1}\Phi_\#-\vk^{-1}R_\#\|_{L^\be}\le 1.
$$
In other words, this reduces the problem of the description of all $p$-best approximants to the problem
of describing all matrix functions $\Q$ in $H^\be(\mm_{n-k.n-k})$ such that 
\bay
\label{NP}
\|\vk^{-1}\Phi_\#-\Q\|_{L^\be}\le 1.
\ey

Note that the problem to describe all $H^\be$ matrix function $\Q$ satisfying \rf{NP} is the classical Nehari problem and in the case of nonuniqueness there is formula parametrizing all solutions.
It was obtained by Adamyan, Arov, and Krein in \cite{AAK1} and \cite{AAK2}
under certain assumptions and by Kheifets
\cite{K} in the most general case; see also Ch. 5 of \cite{P}.

\
\section{\bf $\bs{p}$-superoptimal approximation}
\setcounter{equation}{0}

\

In this section we introduce the notion of $p$-superoptimal approximation and prove that if
$\Phi$ is a rational matrix function then $\Phi$ has a unique $p$-superoptimal approximant.

\medskip

{\bf Definition.} Let $\Phi\in L^p(\mm_{n})\setminus H^p(\mm_{n})$. 
For a function $Q\in H^p(\mm_{n})$, we define the
numbers $\t_j(\Phi,Q)$, $0\le j\le n-1$, by
$$
\t_j(\Phi,Q)=\ess\sup_{\z\in\T}\frac{s_j\big((\Phi-Q)(\z)\big)}{{\rm d}_\Phi(\z)}.
$$
A function $Q\in H^p(\mm_{n})$
is called a $p$-{\it superoptimal approximant} to $\Phi$ if it minimizes lexicographically
the sequence $\t_j(\Phi,Q)$, $0\le j\le\ n-1$.

If $Q$ is a $p$-superoptimal approximant to $\Phi$, we put
$$
\t_j(\Phi)\df\t_j(\Phi,Q),\quad 0\le j\le n-1.
$$

\medskip

Clearly, if $Q$ is a best analytic approximant to $\Phi$ in $L^p(\mm_{m,n})$, then 
$\t_0(\Phi,Q)=1$. It is also clear that if $F$ is a $p$-superoptimal approximant, then
$F$ is a best analytic approximant in $L^p(\mm_{m,n})$.

It is easy to see that if $\Phi$ has gender $k$, then 
$$
\t_j(\Phi)=1\quad\mbox{for}\quad j=0,\cdots,k-1.
$$

In this section we are going to work with rational matrix functions.
When we say that a matrix function defined on the unit circle $\T$ is rational, we mean that it is a restriction of a rational matrix function to the unit circle. It is easy to see that if $A$ is a rational matrix function, then its adjoint $A^*$ is also a rational matrix function.

Suppose now that $\Phi$ is a rational matrix function of size $n\times n$ with no poles in $\T$ and $k$ is the gender of $\Phi$.
As in the proof of Theorem  \ref{ptf}, we consider a co-outer maximizing vector $F$ of the 
Hankel operator $H^{\{k\}}_\Phi$, the matrix function $G$ defined by $G=\bar z\ov{H_\Phi^{\{k\}}F}$, the factorizations
$$
G=G_{\rm co}\cO,\quad F=\U F_{\rm o},\quad\mbox{and}\quad G_{\rm co}=\O G_{\rm o},
$$
where $\cO$ is an inner matrix function of size $k\times k$, $\U$ and $\O$ are inner and co-outer matrix functions
of size $n\times k$, $G_{\rm co}$ is a co-outer matrix functions of size $n\times k$, and
$F_{\rm o}$ and $G_{\rm o}$ are outer matrix functions of size $k\times k$. We also assume that the $k$-balanced matrix functions $\V$ and $\W^{\rm t}$ are given by \rf{cc}, $Q$ is a $p$-best approximant to $\Phi$ and $\Phi-Q$ is factorized as in \rf{pbaf}. Finally, we assume that $F$ is normalized so that
\rf{normal} holds.

\begin{lem}
\label{rat}
Let $\Phi$ be a rational matrix function in $L^p(\mm_n)$. Then the matrix functions $\cO$, $F$, $G$,
$\V$, and $\W$ are also rational.
\end{lem}

\Pf If $\Phi$ is rational, it is easy to see that $H^{\{k\}}_\Phi A$ is rational for an arbitrary 
function $A\in H^q(\bS_2^{n,k})$. In particular, this is true for the function
$H^{\{k\}}_\Phi F$, and so $G$ is rational.

 Let us show that $\cO$ is rational. It is well known (see e.g., \cite{P}, Ch. 2, \S\,5) that  a square inner function $\cU$ is rational if and only if the subspace 
 $$
 K_\cU\df H^2(\C^k)\ominus\cU H^2(\C^k)=\cU H^2_-(\C^k)\bigcap H^2(\C^k)
 $$
 is finite-dimensional. Since $G$ is rational, the Hankel operator 
 $$
 H_{\ov{G}}:H^2(\C^k)\to H^2_-(\C^k)
 $$
 has finite rank (see e.g., \cite{P}, Ch. 2, \S,5). It is easy to see that for $f\in K_{\cO^{\rm t}}$,
 $$
 H_{\ov{G}}f=\ov{G}f.
 $$
 Since $\rank G(\z)=k$ almost everywhere on $\T$, it follows that multiplication
 by $\ov{G}$ has trivial kernel. Thus $K_{\cO^{\rm t}}$ is finite-dimensional, and so $\cO^{\rm t}$ is rational.
 Thus $\cO$ is rational, and so $G_{\rm co}=G\cO^*$ is also rational.

To prove that the matrix function $G_{\rm o}$ is rational, we observe that 
$$
G_{\rm co}^*G_{\rm co}=G_{\rm o}^*G_{\rm o},
$$
and so $G_{\rm o}^*G_{\rm o}$ is a rational matrix function.
The rationality of $G_{\rm o}$ follows now from
the following well-known fact (see \cite{Y}):
if $\Q$ be a matrix outer function of class $H^2(\mm_{k})$, then $\Q$ is rational 
if and only if $\Q^*\Q$ is rational.

 We have $\O=G_{\rm co}G_{\rm o}^{-1}$, and so $\O$
is rational. By Lemma 12.1 of Ch. 14 of \cite{P}, the matrix function $\Xi$ is rational, and so $\W$ is rational.

Let us show that $\V$ is a rational matrix function. Since $\Phi$ is rational, it follows that
$$
\pp_-\Phi=\pp_-(\Phi-Q)=\pp_-\W^*
\left(\begin{matrix}\D&\0\\\0&\Phi_\#\end{matrix}\right)\V^*
$$
is a rational matrix function. Thus
$$
\pp_-\O^{\rm t}\W^*\left(\begin{matrix}\D&\0\\\0&\Phi_\#\end{matrix}\right)\V^*
=\pp_-\O^{\rm t}\pp_-\W^*\left(\begin{matrix}\D&\0\\\0&\Phi_\#\end{matrix}\right)\V^*
$$
is rational. We have
$$
\O^{\rm t}\W^*\left(\begin{matrix}\D&\0\\\0&\Phi_\#\end{matrix}\right)\V^*=
\left(\begin{matrix}\bs{I}_k&\0\end{matrix}\right)\left(\begin{matrix}\D&\0\\\0&\Phi_\#\end{matrix}\right)\V^*
=\D\U^*,
$$
and so $\pp_-\D\U^*$ is a rational matrix function.

Let $h$ be a scalar outer function such that 
$$
|h(\z)|=\left\|G(\z)\right\|_{\bS_2^{n,k}},\quad\z\in\T.
$$
Then $h\in H^2$. It follows from \rf{mv} and \rf{normal} that
\bay
\label{albe}
\|\D(\z)\|_{\mm_{k}}=|h(\z)|^{2/p}\quad\mbox{and}\quad
\left\|F(\z)\right\|_{\bS_2^{n,k}}=|h(\z)|^{2/q}.
\ey
Since $G$ is rational, the function $|h|^2$ is rational.
It follows from the result from \cite{Y} quoted above that the function $h$ is also rational.

Let us show that $F_{\rm o}$ is a maximizing vector of $\bs{H}_\D$ and
$H_\Phi^{\{k\}} F=\ov{\O}\D F_{\rm o}$. Since $F=\U F_{\rm o}$ is a maximizing vector 
of $H^{\{k\}}_\Phi$, we have by Theorem \ref{anal},
\begin{align*}
H^{\{k\}}_\Phi F&=(\Phi-Q)F=\W^*\left(\begin{matrix}\D&\0\\\0&\Phi_\#\end{matrix}\right)\V^*\U F_{\rm o}
\\[.2cm]
&=\left(\begin{matrix}\ov{\O}&\Xi\end{matrix}\right)
\left(\begin{matrix}\D&\0\\\0&\Phi_\#\end{matrix}\right)
\left(\begin{matrix}\U^*\\\Theta^{\rm t}\end{matrix}\right)\U F_{\rm o}
=\ov{\O}\D F_{\rm o}\in H^2_-(\bS_2^{n,k}).
\end{align*}
Since the matrix function $\O$ is co-outer, it follows from Lemma 1.4 of \mbox{Ch. 14} of \cite{P} that 
$\D F_{\rm o}\in H^2_-(\bS_2^{k})$. It is easy to see from \rf{albe} that 
$F_{\rm o}$ is a maximizing vector of $\bs{H}_\D$ and
$H_\Phi^{\{k\}} F=\ov{\O}\D F_{\rm o}$. Thus
$$
\D F_{\rm o}=\bar z\ov{G_{\rm o}}\ov{\cO}.
$$
Consider now the inner-outer factorization of the matrix function $G_{\rm o}\cO$:
$$
G_{\rm o}\cO={\frak O}\cG_{\rm o}.
$$
Clearly, both ${\frak O}$ and $\cG_{\rm o}$ are rational matrix functions. Then
$\D F_{\rm o}=\bar z\ov{{\frak O}\cG_{\rm o}}$, and so
$$
\D=\bar z\ov{{\frak O}\cG_{\rm o}}F_{\rm o}^{-1}
$$
Put 
$$
U=zh^{-2/p}{\frak O}^{\rm t}\D=\ov{\cG}_{\rm o}(h^{2/p}F_{\rm o})^{-1}.
$$
It is easy to see that $U$ is unitary-valued. Put
$$
\Q=h^{2/p}F_{\rm o}\quad\mbox{and}\quad \Q_\#=\cG_{\rm o}^{\rm t}.
$$
Since $U$ is unitary-valued, it is easy to verify that 
$$
\Q^*\Q=\Q_\#\Q_\#^*.
$$
Clearly, $\Q$ is outer.
Since $\Q_\#\Q_\#^*$ is rational, it follows from the result of \cite{Y} quoted above
 that 
$\Q=h^{2/p}F_{\rm o}$ is rational and $U$ is rational. 

We have $\D=\bar z h^{2/p}\ov{\frak O}U$, and so
$$
\pp_-\D\U^*=\pp_-\bar z h^{2/p}\ov{\frak O}U\U^*
=\pp_-\bar z h^{2/p}\ov{\frak O}\Q_\#^*\Q^{-1}\U^*
=\pp_- h^{2/p}\big(\pp_-\bar z\ov{\frak O}\Q_\#^*\Q^{-1}\U^*\big)
$$
is a rational matrix function. Put
$$
{\frak R}\df\pp_-\bar z\ov{\frak O}\Q_\#^*\Q^{-1}\U^*.
$$
Let us show that ${\frak R}$ is rational. We have
$$
\pp_-h{\frak R}=\pp_-h^{2/q}h^{2/p}{\frak R}=\pp_-h^{2/q}\pp_-h^{2/p}{\frak R},
$$
and since $\pp_-h^{2/p}{\frak R}$ is rational and $h^{2/q}\in H^\be$, it follows that  
$\pp_- h{\frak R}$ is rational. Since $h$ is rational and ${\frak R}=\pp_-{\frak R}$, it is easy to see that 
${\frak R}$ is rational.

Finally, since the matrix functions $\bar z\ov{\frak O}\Q_\#^*\Q^{-1}$
and $\pp_-\bar z\ov{\frak O}\Q_\#^*\Q^{-1}\U^*$ are rational, it is easy to verify that
$\U^*$ is rational. Again, it follows from Lemma 12.1 of Ch. 14 of \cite{P} that
$\V$ is rational. $\bl$

To prove the next theorem, we introduce the notation $\L_\a$, $0<\a<1$, for the class of H\"older functions of order $\a$: a function $\f$ on $\T$ is said to belong to the H\"older class $\L_\a$ if
$$
\sup_{\z\ne\t}\frac{|\f(\z)-\f(\t)|}{|\z-\t|^\a}<\be.
$$

In the following theorem we keep all the notation as above.

\begin{thm}
\label{psup}
If $\Phi$ is a rational matrix function, then $\D\in\L_{2/p}(\mm_{k})$
and $h^{-2/p}\Phi_\#\in(H^\be+C)(\mm_{n-k})$.
\end{thm}

\Pf We have 
$$
\D=\bar z h^{2/p}\ov{\frak O}U,
$$
where $\ov{\frak O}U$ is a rational function. If $h$ has no zeros on $\T$, then $\D$ is infinitely differentiable. If $h$ has zeros on $\T$, then, obviously, $h^{2/p}\in \L_{2/p}$, which implies that
$\D\in\L_{2/p}(\mm_{k})$.

Next, since
$$
\W(\Phi-Q)\V=\left(\begin{matrix}\D&\0\\\0&\Phi_\#\end{matrix}\right),
$$
it follows that $\Phi_\#$ is a sum of a rational matrix function and an $H^p$ matrix function. 
Thus there exists a finite Blaschke product $B$ such that $\Phi_\#=\ov{B}\Phi_\heartsuit$,
where \lb$\Phi_\heartsuit\in H^p(\mm_{n-k})$. We also know that 
$\|\Phi_\heartsuit(\z)\|_{\mm_{n-k}}\le|h(\z)|^{2/p}$. Since $h$ is outer, it follows that 
$h^{-2/p}\Phi_\heartsuit\in H^\be(\mm_{n-k})$.
Thus 
$$
h^{-2/p}\Phi_\#=\ov{B}h^{-2/p}\Phi_\heartsuit\in(H^\be+C)(\mm_{n-k}).\quad\bl
$$

\begin{thm}
\label{so}
Let $\Phi$ be an $n\times n$ rational matrix function. Then $\Phi$ has a unique $p$-superoptimal approximant $Q$.
Moreover,
\bay
\label{eq}
\frac{s_j\big((\Phi-Q)(\z)\big)}{{\rm d}_\Phi(\z)}=\t_j(\Phi),\quad 0\le j\le n-1
\ey
almost everywhere on $\T$.
\end{thm}

\Pf Let $R$ be a best analytic approximant to $\Phi$ in $L^p(\mm_{n})$. By Theorem \ref{ptf},
$\Phi-R$ admits a factorization of the form
$$
\Phi-R=\W^*\left(\begin{array}{cc}\D&\0\\\0&\Phi_\#\end{array}\right)\V^*,
$$
where $\V$ and $\W^{\rm t}$ are $k$-balanced matrix function, $k$ is the gender of $\Phi$,
 $\D$ is a \lb$p$-badly approximable
$k\times k$ matrix function function such that the matrix function ${\rm d}_\Phi^{-1}\D$ is unitary-valued,
and $\Phi_\#$ is a matrix function such that 
$$
\|\Phi_\#(\z)\|_{\mm_{n-k}}\le\|\D(\z)\|_{\mm_k},\quad\z\in\T.
$$

It follows from Theorem \ref{paf} that a matrix function \mbox{$Q\in H^p(\mm_{n})$} is a
\lb$p$-superoptimal approximant to $\Phi$ if and only if
$$
\Phi-Q=\W^*\left(\begin{array}{cc}\D&\0\\\0&\Phi_\#-Q_\#\end{array}\right)\V^*,
$$
where $Q_\#\in H^p(\mm_{n-k})$
is a matrix function  such that $h^{-2/p}Q_\#$
is a superoptimal approximant of $h^{-2/p}\Phi_\#$ in $L^\be$. Here $h$ the 
scalar outer function as in 
the proof of Lemma \ref{rat}, i.e., $|h^{2/p}|={\rm d}_\Phi$. By Theorem \ref{psup},
$h^{-2/p}\Phi_\#\in H^\be+C$, and by Theorem 3.3 of Chapter 14 of
\cite{P}, $h^{-2/p}\Phi_\#$ has a unique superoptimal
approximant  in the $L^\be$ norm. 

Formula \rf{eq} is an immediate consequence of Theorem 3.4 of Chapter 14 of \lb \cite{P}.
$\bl$

The following example shows that there are matrix functions in $L^p$, for which there are different
$p$-superoptimal approximants.

\medskip

{\bf Example.} Let $\f$ be a scalar $L^\be$ function such that 
$$
\|\f\|_{L^\be}=\dist_{L^\be}(\f,H^\be)=1,
$$
and such that there is a nonzero best approximant $f\in H^\be$ in the norm of $L^\be$.
It is well known that such functions $\f$ exist (see, e.g., \cite{P}, Ch. 1, \S\,1). Consider the matrix function
$\Phi\in L^\be(\mm_2)$ defined by
$$
\Phi=\left(\begin{matrix}\bar z&\0\\\0&\f\end{matrix}\right).
$$
It is easy to see that both the zero function and the function $\left(\begin{matrix}\0&\0\\\0&f\end{matrix}\right)$
are \lb$p$-superoptimal approximants for any $p\in(2,\be)$.

\

\

\footnotesize
\noindent
\begin{tabular}{p{4.6cm}p{4.5cm}p{4.6cm}}
L. Baratchart  &F.L. Nazarov&V.V. Peller \\
INRIA& Department of Mathematics& Department of Mathematics\\
BP 93& University of Wisconsin&Michigan State University  \\
06902 Sophia-Antipolis Cedex &Madison, Wisconsin 53706& East Lansing, Michigan 48824\\
France&USA&USA
\end{tabular}

\end{document}